\documentclass[12pt]{article}

\usepackage{amstext}
\usepackage{latexsym}
\usepackage{amssymb}

\def\R{\mathbb R}

\def\P{\mathbb P}

\newcommand{\eps}{\varepsilon}

\begin{document}

\newtheorem{theorem}{Theorem}[section]
\renewcommand{\thetheorem}{\arabic{section}.\arabic{theorem}}
\newtheorem{definition}[theorem]{Definition}
\newtheorem{deflem}[theorem]{Definition and Lemma}
\newtheorem{lemma}[theorem]{Lemma}
\newtheorem{example}[theorem]{Example}
\newtheorem{remark}[theorem]{Remark}
\newtheorem{remarks}[theorem]{Remarks}
\newtheorem{cor}[theorem]{Corollary}
\newtheorem{pro}[theorem]{Proposition}
\newtheorem{proposition}[theorem]{Proposition}

\renewcommand{\theequation}{\thesection.\arabic{equation}}

\title{Yet another criterion for global existence 
\\ in the 3D relativistic Vlasov-Maxwell system}

\date{}

\author{{\sc Markus Kunze}\\[2ex]
         Mathematisches Institut, Universit\"at K\"oln, \\
         Weyertal 86-90, D\,-\,50931 K\"oln, Germany \\
         e-mail: mkunze@mi.uni-koeln.de}

\maketitle

\begin{abstract}\noindent
We prove that solutions of the 3D relativistic Vlasov-Maxwell system 
can be extended, as long as the quantity $\sigma_{-1}(t, x)
=\max_{|\omega|=1}\,\int_{\R^3}\frac{dp}{\sqrt{1+p^2}}\,\frac{1}{(1+v\cdot\omega)}\,f(t, x, p)$ 
is bounded in $L^2_x$. 
\end{abstract}


\setcounter{equation}{0}

\section{Introduction and main result}

The relativistic Vlasov-Maxwell system describes the time evolution 
of a plasma, i.e., of an ensemble of charged particles (like ions 
or electrons) in position-momentum phase space $\R^3\times\R^3$. 
Since the particles can move at relativistic speeds the motion of a single particle 
is described by the system 
\begin{equation}\label{charg} 
   \dot{X}=V=\frac{P}{\sqrt{1 + |P|^2}},\quad\dot{V}=E+V\wedge B,
\end{equation} 
where we take only one species for simplicity and choose units where $c=1$ 
for the speed of light. Furthermore, the rest mass and the charge of the particles 
are set to unity. The particle velocity is $V$, whereas $P$ denotes its momentum.  
The vectors $E$ and $B$ in (\ref{charg}) stand for the electric and the magnetic field, respectively. 
Since the number of individual particles in the plasma is large 
one takes a statistical approach and models the time evolution 
by using a density function $f=f(t, x, p)\ge 0$ depending on time $t\in\R$,
position $x\in\R^3$, and momentum $p\in\R^3$. 
Then the requirement that $f$ be constant along the particle trajectories, 
i.e., the solutions of the characteristic equations (\ref{charg}), leads to the Vlasov equation 
\begin{equation}\label{vlas1}
   \partial_t f(t, x, p)+v\cdot\nabla f(t, x, p)+(E(t, x)+v\wedge B(t, x))\cdot\nabla_p f(t, x, p)=0; 
\end{equation}
here $\nabla$ always means $\nabla_x$.
The velocity $v\in\R^3$ associated to $p$ is
\[ v=\frac{p}{\sqrt{1+p^2}}\,,\quad\mbox{thus}\quad p=\frac{v}{\sqrt{1-v^2}}\,, \] 
where $p^2=|p|^2$ and $v^2=|v|^2$ for brevity. The Lorentz force
\[ L=L(t, x, v)=E(t, x)+v\wedge B(t, x)\in\R^3 \]
is obtained from the fields $E$ and $B$, which in turn satisfy
the Maxwell equations
\begin{equation}\label{maxE}
   \partial_t E=\nabla\wedge B-j,\quad\nabla\cdot E=\rho,
\end{equation}
and
\begin{equation}\label{maxB}
   \partial_t B=-\nabla\wedge E,\quad\nabla\cdot B=0.
\end{equation}
The coupling of (\ref{vlas1}) to (\ref{maxE}), (\ref{maxB}) is realized
through the charge density $\rho=\rho(t, x)\in\R$ and the current density
$j=j(t, x)\in\R^3$ via
\begin{equation}\label{rhoj-def}
   \rho(t, x)=\int_{\R^3}f(t, x, p)\,dp\quad\mbox{and}\quad
   j(t, x)=\int_{\R^3} v\,f(t, x, p)\,dp.
\end{equation}
Furthermore, initial data
\begin{equation}\label{inidat}
   f(t=0)=f^{(0)},\quad E(t=0)=E^{(0)},\quad\mbox{and}\quad
   B(t=0)=B^{(0)}
\end{equation}
are prescribed such that the constraint equations
\begin{equation}\label{constr-equ}
   \nabla\cdot E^{(0)}=\rho^{(0)}=\int_{\R^3} f^{(0)}\,dp
   \quad\mbox{and}\quad\nabla\cdot B^{(0)}=0
\end{equation}
are satisfied. Good general introductions to the subject can be found in \cite{strauss, glassey}. 
\medskip 

The relativistic Vlasov-Maxwell system comprises a complicated system of nonlinear partial differential equations. 
Local existence of solutions for smooth and compactly supported (or sufficiently decaying) initial data 
and a sufficient condition for global existence has been know for some time. 
More precisely, we have the following result. 

\begin{theorem}[Glassey/Strauss]\label{locexi} Let the initial data
\[ f^{(0)}\in C_0^1(\R^3\times\R^3)\quad\mbox{and}\quad
   E^{(0)}, B^{(0)}\in C^2_b(\R^3; \R^3)\cap L^2(\R^3; \R^3) \]
be given such that the constraint equations (\ref{constr-equ}) are satisfied.
Then there exists a maximal local solution $(f, E, B)$ on a time interval $[0, T_{{\rm max}}[$
to the relativistic Vlasov-Maxwell system (\ref{vlas1}), (\ref{maxE}), (\ref{maxB}),
(\ref{rhoj-def}) such that the initial data are attained at $t=0$; see (\ref{inidat}).
Furthermore, if there is a function $\varpi\in C([0, \infty[)$ so that
\begin{equation}\label{Pinfty-crit}
   P_\infty(t)\le\varpi(t),\quad t\in [0, T_{{\rm max}}[,
\end{equation}
then $T_{{\rm max}}=\infty$.
\end{theorem}
Here
\[ P_\infty(t)=10+\sup\,\Big\{|p|:\,\,\exists\,s\in [0, t]\,\,\exists\,x\in\R^3:
   \,\,f(s, x, p)\neq 0\Big\} \]
is the maximal momentum up to the time $t$. For a proof, see \cite{glstr}, 
or \cite{klasta, bogopa} for the same result obtained by different methods. 
The work \cite{glstr} has been generalized to initial data of non-compact support 
in \cite{glstrtokyo, kunze}, the latter is also extending \cite{horst, rein}. 
The problem of unrestricted global existence has been studied by many people. 
Only in the framework of weak solution it has been solved in \cite{diplions}. 
Regarding classical solutions, adding spherical symmetry \cite{horst}, 
smallness of initial data \cite{glstrabs}, or ``near neutrality'' \cite{GlSchnn, rein} 
has turned out to be sufficient to close the case. 
Another remarkable work is \cite{glschaeff}, where global existence 
was shown for the ``two and one-half-dimensional'' system, i.e., $x\in\R^2$ and $p\in\R^3$. 
\medskip 

Due to the lack of a general global result in 3D it is natural to focus on deriving further 
continuation criteria, apart from (\ref{Pinfty-crit}), which might be easier to check 
(although, of course, all criteria would be equivalent in the end if global existence was known). 
The following result summarizes some early attempts in this direction. 

\begin{theorem} The following are equivalent: 
\begin{itemize}
\item[(a)] $T_{{\rm max}}=\infty$. 
\item[(b)] There is a function $\varpi_1\in C([0, \infty[)$ such that 
$P_\infty(t)\le\varpi_1(t)$ for $t\in [0, T_{{\rm max}}[$. 
\item[(c)] There is a function $\varpi_2\in C([0, \infty[)$ such that
\[ \sup\bigg\{\int_{\R^3}\sqrt{1+p^2}\,f(t, x, p)\,dp: x\in\R^3\bigg\}\le\varpi_2(t),
   \quad t\in [0, T_{{\rm max}}[. \] 
\item[(d)] There is a function $\varpi_3\in C([0, \infty[)$ such that
\[ \sup\bigg\{|P(t; 0, x, p)-p|: (x, p)\in\R^3\times\R^3\bigg\}\le\varpi_3(t),
   \quad t\in [0, T_{{\rm max}}[. \] 
Here $s\mapsto (X(s; t, x, p), P(s; t, x, p))$ is the solution to the characteristic system 
(\ref{charg}) which at $s=t$ equals $(x, p)$. 
\item[(e)] 
\begin{itemize}
\item[(i)] $f^{(0)}=0$ or 
\item[(ii)] there is $\eps>0$ and $R_0>0$ such that $\int_{|x|\le R_0}\int_{\R^3} f^{(0)}(x, p)\,dx\,dp>\eps$ 
so that for all $R\in [0, R_0+T_{{\rm max}}[$ it holds that 
\[ \lim_{t\to T_{{\rm max}}}\int_{|x|\le R}\bigg[\int_{\R^3}\sqrt{1+p^2}\,f(t, x, p)\,dp
   +\frac{1}{2}(|E(t, x)|^2+|B(t, x)|^2)\bigg]\,dx=0. \] 
\end{itemize}   
\end{itemize}
\end{theorem} 
Parts (c) and (e) are more or less contained in \cite{glstrtokyo} and \cite{glstrcontemp}, respectively; 
see \cite{kunze} for a detailed proof. After some dormant period recently the interest in the subject 
has been revived and further criteria have been obtained, using refined techniques. 
To state the results we need to introduce the quantity 
\begin{equation}\label{calidef} 
   {\cal I}_\theta(t, x)=\int_{\R^3} (1+p^2)^{\frac{\theta}{2}} f(t, x, p)\,dp
\end{equation} 
for $\theta>0$. 

\begin{theorem}\label{criti2} 
The following are equivalent: 
\begin{itemize}
\item[(a)] $T_{{\rm max}}=\infty$. 
\item[(f)] There is a function $\varpi_4\in C([0, \infty[)$ such that 
\[ \sup\,\Big\{|\P_Q p|:\,\,\exists\,s\in [0, t]\,\,\exists\,x\in\R^3:
   \,\,f(s, x, p)\neq 0\Big\}\le\varpi_4(t),
   \quad t\in [0, T_{{\rm max}}[. \] 
Here $Q\subset\R^3$ is a two-dimensional plane in $p$-space containing the origin 
and $\P$ denotes the projection onto $Q$. 
\item[(g)] There is a function $\varpi_5\in C([0, \infty[)$ such that 
\[ {\|\rho(t)\|}_{L_x^\infty(\R^3)}={\|{\cal I}_0(t)\|}_{L_x^\infty(\R^3)}\le\varpi_5(t),
   \quad t\in [0, T_{{\rm max}}[.  \] 
\item[(h)] Let $q\in ]2, \infty]$ and $\theta>2/q$, or $q\in [1, 2]$ and $\theta>8/q-3$. 
Then there is a function $\varpi_6\in C([0, \infty[)$ such that 
\[ {\|{\cal I}_\theta(t)\|}_{L_x^q(\R^3)}\le\varpi_6(t),
   \quad t\in [0, T_{{\rm max}}[.  \]
\end{itemize}
\end{theorem} 
Part (f), in fact for more general data, is due to \cite{LStrgd}. 
It shows that one does not have to control the full momentum support, 
but only its projection to some plane through the origin. Part (g) has been proved in \cite{sosill}. 
The most recent result is (h), which is cited from \cite{LStrbd}, and once again 
holds for more general data. The results from \cite{LStrbd} generalize those of 
\cite{pallard}, where $q\in [6, \infty[$ and $\theta>4/q$, or $q<6$ and $\theta>22/q-3$ was assumed. 
\medskip 

In this work we propose to study another quantity, which is 
\begin{equation}\label{sig-1-def} 
   \sigma_{-1}(t, x)=\max_{|\omega|=1}\,\int_{\R^3}\frac{dp}{\sqrt{1+p^2}}\,\frac{1}{(1+v\cdot\omega)}\,f(t, x, p)
\end{equation} 
and which comes up naturally in the course of the estimates. Our main result is as follows. 

\begin{theorem}\label{sig-1thm} Suppose that the initial data
\[ f^{(0)}\in C_0^1(\R^3\times\R^3)\quad\mbox{and}\quad
   E^{(0)}, B^{(0)}\in C^2_b(\R^3; \R^3)\cap L^2(\R^3; \R^3) \]
are given such that the constraint equations (\ref{constr-equ}) are satisfied.
Let $\sigma_{-1}$ be defined by (\ref{sig-1-def}).
If there is a function $\varpi\in C([0, \infty[)$ so that
\begin{equation}\label{criter-hyp}
   {\|\sigma_{-1}\|}_{L^\infty_t L^2_x(S_T)}\le\varpi(T),
   \quad T\in [0, T_{{\rm max}}[,
\end{equation}
for $S_T=[0, T]\times\R^3$, then $T_{{\rm max}}=\infty$.
\end{theorem}

The point is that $\int_{\R^3}\frac{dp}{\sqrt{1+p^2}}\,f\in L^\infty_t L^2_x$ 
by energy conservation, so as compared to $\sigma_{-1}$ one might hope to 
``get away with a logarithmic loss''. The method of proof is somewhat similar to \cite{LStrbd}, 
in that Strichartz estimates for a wave equation related to $E$ and $B$ are applied. 
However, we can avoid the use of iteration sequences and bounds on the field derivatives, 
which makes the argument more direct. Comparing ${\|\sigma_{-1}(t)\|}_{L^2_x(\R^3)}$ 
to ${\|{\cal I}_\theta(t)\|}_{L_x^q(\R^3)}$ one can also derive a corollary 
in the fashion of Theorem \ref{criti2}(h). 
 
\begin{cor}\label{thmcor} Under the hypotheses of Theorem \ref{sig-1thm}, let $\theta>1$ 
and $q\in ]\frac{4}{\theta+1}, \infty[$ be given. Then $T_{{\rm max}}=\infty$ is equivalent 
to the existence of a function $\varpi_7\in C([0, \infty[)$ such that 
\[ {\|{\cal I}_\theta(t)\|}_{L_x^q(\R^3)}\le\varpi_7(t),
   \quad t\in [0, T_{{\rm max}}[.  \]
In particular, this gives something new as compared to Theorem \ref{criti2}(h) 
for $\theta>1$ and $q\in ]\frac{4}{\theta+1}, \frac{8}{\theta+3}]$. 
\end{cor} 

The paper is organized as follows. In Section \ref{prel-sect} we collect some preliminary results 
which are well-known in general. Then we turn to deriving suitable bounds on $E$ and $B$ 
in Section \ref{EB-bound}; they mainly rely on the representation formulae for the fields 
due to Glassey and Strauss and on Strichartz estimates for the wave equation. 
The argument for the proof of Theorem \ref{sig-1thm} is elaborated in Section \ref{mainprf-sect}. 
Finally, Section \ref{prcor-sect} contains the proof of Corollary \ref{thmcor}. 
\smallskip

Constants which only depend on the initial data are denoted by $C(0)$, 
whereas $C$ is a numerical constant. Sometimes data terms are not made explicit 
and are only written as ``$({\rm data})$''. By our hypotheses they are good enough 
with regard to all the estimates we will be aiming for. 
\smallskip

The wave operator on $\R\times\R^3$ is $\square=\partial_t^2-\Delta$.
For functions $h=h(t, x)$ define
\begin{eqnarray}\label{waveop-def}
   (\square^{-1}h)(t, x) & = & \int_0^t\frac{ds}{4\pi s}
   \,\int_{|y-x|=s} dS(y)\,h(t-s, y)
   =\int_{|y-x|\le t}\frac{dy}{4\pi |x-y|}\,h(t-|x-y|, y)
   \nonumber\\ & = & \int_{|y|\le t}\frac{dy}{4\pi |y|}\,h(t-|y|, x+y).
\end{eqnarray}
Then $g=\square^{-1}h$ is the unique solution to
\[ \square\,g=h,\quad g(0)=\partial_t g(0)=0. \]


\setcounter{equation}{0}

\section{Some preliminaries}
\label{prel-sect} 

From the system (\ref{vlas1}), (\ref{maxE}), and (\ref{maxB}) it follows that the energy 
\[ {\cal E}(t)=\int_{\R^3}\int_{\R^3}\sqrt{1+p^2}\,f(t, x, p)\,dx\,dp
   +\frac{1}{2}\int_{\R^3}(|E(t, x)|^2+|B(t, x)|^2)\,dx \] 
is conserved along sufficiently regular solutions (as we are dealing with); 
note that $\nabla_p\sqrt{1+p^2}=v$. In addition, (\ref{rhoj-def}) and (\ref{vlas1})
yield the continuity equation $\partial_t\rho+\nabla\cdot j=0$. 

For $(t, x, p)$ fixed let $s\mapsto (X(s; t, x, p), P(s; t, x, p))$
denote the solution of the characteristic initial value problem (\ref{charg}) 
which at $s=t$ equals $(x, p)$. Then (\ref{vlas1}) is equivalent to
\[ \frac{d}{ds}\,f(s, X(s, t, x, p), P(s, t, x, p))=0, \]
which leads to the relation
\begin{equation}\label{ff0}
   f(t, x, p)=f^{(0)}(X(0, t, x, p), P(0, t, x, p))
\end{equation}
for the initial data $f^{(0)}(x, p)=f(0, x, p)$. Thus in particular 
\[ {\cal L}(t)={\|f(t)\|}_{L^\infty(\R^3\times\R^3)}
   ={\|f^{(0)}\|}_{L^\infty(\R^3\times\R^3)}={\cal L}(0). \] 
Since every map $(x, p)\mapsto (X(s; t, x, p), P(s; t, x, p))$ 
is a measure preserving diffeomorphism of $\R^3\times\R^3$, it follows from (\ref{ff0}) 
that for instance 
\[ {\|\rho(t)\|}_{L^1(\R^3)}=\int_{\R^3}\int_{\R^3}f(t, x, p)\,dx\,dp
   =\int_{\R^3}\int_{\R^3}f^{(0)}(\bar{x}, \bar{p})\,d\bar{x}\,d\bar{p}
   ={\|\rho^{(0)}\|}_{L^1(\R^3)}, \] 
where $\rho^{(0)}(x)=\int_{\R^3}f^{(0)}(x, p)\,dp$, which expresses the conservation of mass.
\medskip 

The following result is also known, but we nevertheless include a proof 
in order to make the presentation self-contained.  

\begin{lemma}\label{moment-ode} Let $m_k$ be defined by 
\[ m_k(t)=1+\int_{\R^3}\int_{\R^3} (1+p^2)^\frac{k}{2} f(t, x, p)\,dx\,dp. \]
If $k\in [2, \infty[$, then
\[ m_k(t)\le m_k(0)+C{\cal L}(0)^{\frac{1}{k+3}}\,\int_0^t {\|E(s)\|}_{L^{k+3}_x(\R^3)}
   \,m_k(s)^{\frac{k+2}{k+3}}\,ds, \]
where $C$ depends on $k$.
\end{lemma}
{\bf Proof\,:} The Vlasov equation (\ref{vlas1}) yields
\begin{eqnarray*}
   \frac{dm_k}{dt} & = & \int_{\R^3}\int_{\R^3} (1+p^2)^\frac{k}{2}\,\partial_t f\,dx\,dp
   =-\int_{\R^3}\int_{\R^3} (1+p^2)^\frac{k}{2}\,\nabla_p\cdot((E+v\wedge B)f)\,dx\,dp
   \\ & = & k\int_{\R^3}\int_{\R^3} (1+p^2)^{\frac{k}{2}-1}\,p\cdot ((E+v\wedge B)f)\,dx\,dp
   \\ & = & k\int_{\R^3}\int_{\R^3} (1+p^2)^\frac{k-1}{2}\,(v\cdot E)f\,dx\,dp,
\end{eqnarray*}
so that
\begin{equation}\label{rutga}
   \Big|\frac{dm_k}{dt}\Big|
   \le k\,{\|E\|}_{L^{q'}_x(\R^3)}
   {\bigg\|\int_{\R^3} (1+p^2)^\frac{k-1}{2}\,f\,dp\,\bigg\|}_{L^q_x(\R^3)}
\end{equation}
for $q\in [1, \infty]$. If $R\in ]0, \infty[$ and $\theta\in ]0, \infty[$ are fixed, then
\begin{eqnarray*}
   \int_{\R^3} (1+p^2)^\frac{k-1}{2}\,f\,dp
   & = & \int_{|p|\le R} (1+p^2)^\frac{k-1}{2}\,f\,dp
   +\int_{|p|>R} (1+p^2)^\frac{k-1}{2}\,f\,dp
   \\ & \le & {\cal L}(0)\int_{|p|\le R} (1+p^2)^\frac{k-1}{2}\,dp
   +R^{-2\theta}\int_{|p|>R} (1+p^2)^{\frac{k-1}{2}+\theta}\,f\,dp
   \\ & \le & C{\cal L}(0)(1+R)^{k+2}
   +R^{-2\theta}\int_{\R^3} (1+p^2)^{\frac{k-1}{2}+\theta}\,f\,dp,
\end{eqnarray*}
where $C$ depends on $k$. Tacitly assuming $R\in [1, \infty[$ for the optimal $R$
(otherwise the compact $x$-support is useful to obtain the needed bound), this yields
\[ \int_{\R^3} (1+p^2)^\frac{k-1}{2}\,f\,dp\le C{\cal L}(0)^{\frac{2\theta}{k+2(1+\theta)}}
   \,\bigg(\int_{\R^3} (1+p^2)^{\frac{k-1}{2}+\theta}\,f\,dp\bigg)^{\frac{k+2}{k+2(1+\theta)}}, \]
where $C$ depends on $k$ and $\theta$. Therefore the estimate
\[ {\bigg\|\int_{\R^3} (1+p^2)^\frac{k-1}{2}\,f\,dp\,\bigg\|}_{L^{\frac{k+2(1+\theta)}{k+2}}_x(\R^3)}
   \le C{\cal L}(0)^{\frac{2\theta}{k+2(1+\theta)}}\,m_{k-1+2\theta}^{\frac{k+2}{k+2(1+\theta)}} \]
is found. Putting $q=1+\frac{2\theta}{k+2}$, from $k-1+2\theta=(k+2)q-3$ it follows that
\[ {\bigg\|\int_{\R^3} (1+p^2)^\frac{k-1}{2}\,f\,dp\,\bigg\|}_{L^q_x(\R^3)}
   \le C{\cal L}(0)^{\frac{1}{q'}}\,m_{(k+2)q-3}^{\frac{1}{q}} \] 	
is verified for $k\in [2, \infty[$ and $q\in ]1, \infty[$, where $C$ depends on $k$ and $q$.
Using this in (\ref{rutga}),
\[ \Big|\frac{dm_k}{dt}\Big|
   \le C{\cal L}(0)^{\frac{1}{q'}}\,{\|E\|}_{L^{q'}_x(\R^3)}
   \,m_{(k+2)q-3}^{\frac{1}{q}}. \]
For the particular choice of $q=\frac{k+3}{k+2}$ and $q'=k+3$, the claimed bound is obtained.
{\hfill$\Box$}\bigskip

\begin{lemma}\label{sig-1einP}
If $q\in ]1, \infty[$ and $\alpha\in [0, 1-\frac{1}{2q}[$, then
\[ {\|\sigma_{-1}(t)\|}_{L^q_x(\R^3)}\le CP_\infty(t)^{2(1-\alpha)-\frac{1}{q}}
   \,m_{2\alpha q}(t)^{\frac{1}{q}}, \]
where $C$ depends on $q$, $\alpha$, and ${\cal L}(0)$.
\end{lemma}
{\bf Proof\,:} Fix $\omega\in\R^3$ such that $|\omega|=1$.
Then by H\"older's inequality and by Lemma \ref{efu-lem1}(b) below 
for $\theta=(\alpha+\frac{1}{2})q'$ and $\kappa=q'$,
\begin{eqnarray*}
   \lefteqn{\int_{\R^3}\frac{dp}{\sqrt{1+p^2}}\,\frac{1}{(1+v\cdot\omega)}\,f(t)}
   \\ & = & \int_{\R^3}\frac{dp}{(1+p^2)^{\alpha+\frac{1}{2}}}\,\frac{1}{(1+v\cdot\omega)}
   \,(1+p^2)^\alpha f(t)
   \\ & \le & \bigg(\int_{|p|\le P_\infty(t)}
   \frac{dp}{(1+p^2)^{(\alpha+\frac{1}{2})q'}}\,\frac{1}{(1+v\cdot\omega)^{q'}}\bigg)^{1/q'}
   \bigg(\int_{\R^3}\,(1+p^2)^{\alpha q} f(t)^q\,dp\bigg)^{1/q}
   \\ & \le & C(\alpha, q)\,{\cal L}(0)^{\frac{q-1}{q}}\,P_\infty(t)^{2(1-\alpha)-\frac{1}{q}}
   \,\bigg(\int_{\R^3}\,(1+p^2)^{\alpha q} f(t)\,dp\bigg)^{1/q}.
\end{eqnarray*}
Taking first the $\max_{|\omega|=1}$, then the $q$'th power,
and then integrating $\int_{\R^3} dx$, the claimed bound
is obtained. {\hfill$\Box$}\bigskip

The following general integration lemma is useful at many places. 

\begin{lemma}\label{efu-lem1} Suppose that $R\in [10, \infty[$ and $|\omega|=1$.
\begin{itemize}
\item[(a)] If $\theta\in [0, \frac{3}{2}[$, then
\[ \int_{|p|\le R}\frac{dp}{(1+p^2)^\theta}\,\frac{1}{1+v\cdot\omega}
   \le C\ln R\,R^{3-2\theta}, \]
where $C$ depends on $\theta$.
\item[(b)] If $\theta, \kappa\in [0, \infty[$ are such that
$\theta<\kappa+\frac{1}{2}$ and $\kappa>1$, then
\[ \int_{|p|\le R}\frac{dp}{(1+p^2)^\theta}\,\frac{1}{(1+v\cdot\omega)^\kappa}
   \le CR^{1+2(\kappa-\theta)}, \]
where $C$ depends on $\theta$ and $\kappa$.
\item[(c)] If $\theta\in [0, 1[$ and $|v|<1$, then
\[ \int_{|\omega|=1}\frac{dS(\omega)}{(1+v\cdot\omega)^\theta}\le C, \]
where $C$ depends on $\theta$.
\end{itemize}
\end{lemma}
{\bf Proof\,:} (a) First $\omega$ is rotated to $(0, 0, 1)$.
Then spherical coordinates and the transformation 
\begin{equation}\label{sigma-r}
   \sigma=\frac{r}{\sqrt{1+r^2}},\quad r=\frac{\sigma}{\sqrt{1-\sigma^2}},
   \quad d\sigma=(1-\sigma^2)^{3/2} dr,\quad 1+r^2=(1-\sigma^2)^{-1},
\end{equation} 
are used to get
\begin{eqnarray*}
   \int_{|p|\le R}\frac{dp}{(1+p^2)^\theta}\,\frac{1}{1+v\cdot\omega}
   & = & \int_{|p|\le R}\frac{dp}{(1+p^2)^\theta}\,\frac{1}{1+v_3}
   \le C\int_0^R\frac{dr\,r^2}{(1+r^2)^\theta}\,\frac{1}{1+\frac{r\cos\varphi}{\sqrt{1+r^2}}}
   \\ & = & C\int_0^{R^\flat}\frac{d\sigma\,\sigma^2}{(1-\sigma^2)^{\frac{5}{2}-\theta}}\,
   \int_{-1}^1\frac{ds}{1+\sigma s}
   \\ & = & C\int_0^{R^\flat}\frac{d\sigma\,\sigma}{(1-\sigma^2)^{\frac{5}{2}-\theta}}\,
   \ln\Big(\frac{1+\sigma}{1-\sigma}\Big).
\end{eqnarray*}
Since $\ln(\frac{1+\sigma}{1-\sigma})\le\ln(\frac{4}{1-\sigma^2})
\le\ln(\frac{4}{1-(R^\flat)^2})=\ln(4(1+R^2))\le\ln(8R^2)\le 3\ln R$,
the claim follows. (b) Similar as in (a),
\begin{eqnarray*}
   \int_{|p|\le R}\frac{dp}{(1+p^2)^\theta}\,\frac{1}{(1+v\cdot\omega)^\kappa}
   & \le & C\int_0^{R^\flat}\frac{d\sigma\,\sigma^2}{(1-\sigma^2)^{\frac{5}{2}-\theta}}\,
   \int_{-1}^1\frac{ds}{(1+\sigma s)^\kappa}
   \\ & \le & C\int_0^{R^\flat}\frac{d\sigma\,\sigma}{(1-\sigma^2)^{\frac{5}{2}-\theta}}\,
   \frac{1}{(1-\sigma)^{\kappa-1}}
   \\ & \le & C\int_0^{R^\flat}\frac{d\sigma\,\sigma}{(1-\sigma^2)^{\frac{3}{2}-\theta+\kappa}}
   \le CR^{1+2(\kappa-\theta)}.
\end{eqnarray*}
(c) First consider the case where $|v|\le 1/2$. Then $1+v\cdot\omega\ge 1-|v|\ge 1/2$ yields
\[ \int_{|\omega|=1}\frac{dS(\omega)}{(1+v\cdot\omega)^\theta}\le 4\pi\,2^\theta. \]
If $|v|\ge 1/2$, then $v$ is rotated to $(0, 0, |v|)$ to get
\begin{eqnarray*}
   \int_{|\omega|=1}\frac{dS(\omega)}{(1+v\cdot\omega)^\theta}
   & = & \int_{|\omega|=1}\frac{dS(\omega)}{(1+|v|\,\omega_3)^\theta}
   =\int_0^{2\pi} d\varphi\int_0^\pi d\theta
   \,\frac{\sin\theta}{(1+\cos\theta\,|v|)^\theta}
   \\ & = & 2\pi\int_{-1}^1\,\frac{ds}{(1+s\,|v|)^\theta}
   =\frac{2\pi}{(1-\theta)|v|}\,\Big((1+|v|)^{1-\theta}-(1-|v|)^{1-\theta}\Big)
   \\ & \le & \frac{4\pi\,2^{1-\theta}}{(1-\theta)}.
\end{eqnarray*}
This completes the proof. {\hfill$\Box$}\bigskip

 
\section{Bounds on the fields}
\label{EB-bound}

First we recall the following representation of the fields $E$ and $B$
from \cite[(A13), (A14), (A3)]{schaeffer:86}.
\begin{eqnarray}
   E & = & E_D+E_{DT}+E_\flat+E_\sharp,
   \label{E-form}
   \\ B & = & B_D+B_{DT}+B_\flat+B_\sharp,
   \label{B-form}
\end{eqnarray}
where
\begin{eqnarray}
   E_D(t, x) & = &
   \partial_t\bigg(\frac{t}{4\pi}\int_{|\omega|=1} E^{(0)}(x+t\omega)\,d\omega\bigg)
   +\frac{t}{4\pi}\int_{|\omega|=1}\partial_t E(0, x+t\omega)\,d\omega
   \quad {\rm (data)}, 
   \nonumber \\
   E_{DT}(t, x) & = & -\frac{1}{t}\int_{|y|=t}\int_{\R^3} K_{E,\,DT}(\omega, v)
   f^{(0)}(x+y, p)\,dp\,d\sigma(y)
   \quad {\rm (data)}, 
   \nonumber \\
   E_\flat(t, x) & = & -\int_{|y|\le t}\frac{dy}{|y|^2}
   \int_{\R^3} dp\,K_{E,\,\flat}(\omega, v) f(t-|y|, x+y, p), 
   \label{Eflat-def} \\
   E_\sharp(t, x) & = & -\int_{|y|\le t}\frac{dy}{|y|}
   \int_{\R^3} dp\,K_{E,\,\sharp}(\omega, v)\,(Lf)(t-|y|, x+y, p),
   \nonumber
\end{eqnarray}
and
\begin{eqnarray*}
   B_D(t, x) & = &
   \partial_t\bigg(\frac{t}{4\pi}\int_{|\omega|=1} B^{(0)}(x+t\omega)\,d\omega\bigg)
   +\frac{t}{4\pi}\int_{|\omega|=1}\partial_t B(0, x+t\omega)\,d\omega
   \quad {\rm (data)}, 
   \\ B_{DT}(t, x) & = & \frac{1}{t}\int_{|y|=t}\int_{\R^3} K_{B,\,DT}(\omega, v)
   f^{(0)}(x+y, p)\,dp\,d\sigma(y)
   \quad {\rm (data)}, 
   \\ B_\flat(t, x) & = & \int_{|y|\le t}\frac{dy}{|y|^2}
   \int_{\R^3} dp\,K_{B,\,\flat}(\omega, v) f(t-|y|, x+y, p), 
   \\ B_\sharp(t, x) & = & \int_{|y|\le t}\frac{dy}{|y|}
   \int_{\R^3} dp\,K_{B,\,\sharp}(\omega, v)\,(Lf)(t-|y|, x+y, p),
\end{eqnarray*}
defining $\omega=|y|^{-1}y$. The respective kernels are given by
\begin{eqnarray*}
   K_{E,\,DT}(\omega, v) & = & (1+v\cdot\omega)^{-1}
   (\omega-(v\cdot\omega)v), \\
   K_{E,\,\flat}(\omega, v)  & = & (1+p^2)^{-1}(1+v\cdot\omega)^{-2}(v+\omega), \\
   K_{E,\,\sharp}(\omega, v)  & = & (1+p^2)^{-1/2}(1+v\cdot\omega)^{-2}
   \Big[(1+v\cdot\omega)+((v\cdot\omega)\omega-v)\otimes v
   \\ & & \hspace{14em} -(v+\omega)\otimes\omega\Big]\in\R^{3\times 3},
\end{eqnarray*}
and
\begin{eqnarray*}
   K_{B,\,DT}(\omega, v) & = & -(1+v\cdot\omega)^{-1}(v\wedge\omega), \\
   K_{B,\,\flat}(\omega, v) & = & -(1+p^2)^{-1}(1+v\cdot\omega)^{-2}(v\wedge\omega), \\
   K_{B,\,\sharp}(\omega, v) & = & (1+p^2)^{-1/2}(1+v\cdot\omega)^{-2}
   \Big[(1+v\cdot\omega)\,\omega\wedge (\ldots)
   \\ & & \hspace{14em} -(v\wedge\omega)\otimes (v+\omega)\Big]\in\R^{3\times 3}.
\end{eqnarray*}

\begin{lemma}\label{kern-bd} The following (known) estimates hold.
\begin{eqnarray}
   & & |K_{E,\,DT}(\omega, v)|+|K_{B,\,DT}(\omega, v)|\le C(1+v\cdot\omega)^{-1/2},
   \nonumber\\[1ex] & & |K_{E,\,\flat}(\omega, v)|+|K_{B,\,\flat}(\omega, v)|
   \le C(1+p^2)^{-1}(1+v\cdot\omega)^{-3/2},
   \nonumber\\[1ex] & & |K_{E,\,\sharp}(\omega, v)z|+|K_{B,\,\sharp}(\omega, v)z|
   \le C(1+p^2)^{-1/2}(1+v\cdot\omega)^{-1}|z|
   \quad (z\in\R^3).
   \label{anglr}
\end{eqnarray}
\end{lemma}
{\bf Proof\,:} The first two lines are a consequence of
\begin{eqnarray}
   |\omega-(v\cdot\omega)v| & = & \Big(1-2(v\cdot\omega)^2+(v\cdot\omega)^2 v^2\Big)^{1/2}
   \le\Big(1-(v\cdot\omega)^2\Big)^{1/2}\le\sqrt{2}\,(1+v\cdot\omega)^{1/2},
   \nonumber
   \\ |v+\omega| & = & \Big(v^2+2(v\cdot\omega)+1\Big)^{1/2}\le\sqrt{2}\,(1+v\cdot\omega)^{1/2},
   \label{vecnorm1}
   \\ |v\wedge\omega| & = & |(v+\omega)\wedge\omega|\le |v+\omega|\le\sqrt{2}\,(1+v\cdot\omega)^{1/2}.
   \nonumber
\end{eqnarray}
The bound on $|K_{B,\,\sharp}(\omega, v)z|$ is immediate from the preceding estimates.
To bound $|K_{E,\,\sharp}(\omega, v)z|$, finally note that
\begin{eqnarray*}
   \Big[((v\cdot\omega)\omega-v)\otimes v-(v+\omega)\otimes\omega\Big]z
   & = & (v\cdot z)((v\cdot\omega)\omega-v)-(\omega\cdot z)(v+\omega)
   \\ & = & -\,(\omega-(v\cdot\omega)v)\cdot z\,(v+\omega)
   -(1+v\cdot\omega)(v\cdot z)\,v.
\end{eqnarray*}
This yields the claim.
{\hfill$\Box$}\bigskip

For functions $h=h(t, x)$ define the operator ${\cal W}$ by
\begin{eqnarray}\label{calW-def}
   ({\cal W}h)(t, x) & = & \int_0^t\frac{ds}{4\pi s^2}
   \,\int_{|y-x|=s} dS(y)\,h(t-s, y)
   =\int_{|y-x|\le t}\frac{dy}{4\pi |x-y|^2}\,h(t-|x-y|, y)
   \nonumber\\ & = & \int_{|y|\le t}\frac{dy}{4\pi |y|^2}\,h(t-|y|, x+y).
\end{eqnarray}

\begin{lemma} The following estimates hold.
\begin{eqnarray}
   & & |E_D(t, x)|+|E_{DT}(t, x)|+|B_D(t, x)|+|B_{DT}(t, x)|\le C\,{\rm (data)},
   \label{EB-DDT-esti}
   \\[1ex] & & |E_\flat(t, x)|+|B_\flat(t, x)|\le C\,({\cal W}\sigma_{-1})(t, x),
   \label{EB-T-esti}
   \\[1ex] & & |E_\sharp(t, x)|+|B_\sharp(t, x)|
   \le C\,\Big(\square^{-1}((|E|+|B|)\sigma_{-1})\Big)(t, x).
   \label{EB-S-esti}
\end{eqnarray}
\end{lemma}
{\bf Proof\,:} Concerning the second pair of estimates, by Lemma \ref{kern-bd} for instance
\begin{eqnarray*}
   |E_\flat(t, x)| & \le & C\int_{|y|\le t}\frac{dy}{|y|^2}
   \int\frac{dp}{1+p^2}\,\frac{1}{(1+v\cdot\omega)^{3/2}}\,f(t-|y|, x+y, p)
   \\ & \le & C\int_{|y|\le t}\frac{dy}{|y|^2}
   \int\frac{dp}{\sqrt{1+p^2}}\,\frac{1}{(1+v\cdot\omega)}\,f(t-|y|, x+y, p)
   \\ & \le & C\int_{|y|\le t}\frac{dy}{|y|^2}\,\sigma_{-1}(t-|y|, x+y)
   =C\,({\cal W}\sigma_{-1})(t, x),
\end{eqnarray*}
using the trivial bound $1+v\cdot\omega
\ge 1-|v|\ge\frac{1}{2}(1-v^2)=\frac{1}{2(1+p^2)}$,
so that $(1+v\cdot\omega)^{-1/2}\le\sqrt{2}\,\sqrt{1+p^2}$. The same argument
can be used to show that also $|B_\flat(t, x)|\le C({\cal W}\sigma_{-1})(t, x)$.
For $E_\sharp$, again Lemma \ref{kern-bd} may be invoked to give
\begin{eqnarray*}
   |E_\sharp(t, x)| & \le & C\int_{|y|\le t}\frac{dy}{|y|}
   \int\frac{dp}{\sqrt{1+p^2}}\,\frac{1}{(1+v\cdot\omega)}\,(|L|f)(t-|y|, x+y, p)
   \\ & \le & C\int_{|y|\le t}\frac{dy}{|y|}\,(|E|+|B|)(t-|y|, x+y)
   \,\sigma_{-1}(t-|y|, x+y)
   \\ & = & C\,\Big(\square^{-1}((|E|+|B|)\sigma_{-1})\Big)(t, x),
\end{eqnarray*}
recall (\ref{waveop-def}). The bound on $|B_\sharp(t, x)|$ is analogous.
{\hfill$\Box$}\bigskip

For the wave equation the following Strichartz estimates are known;
see \cite[(4.9), p.~100]{sogge}. For every $\gamma\in ]0, 1[$
there is a constant $C^\ast_\gamma>0$ with the following property.
Let $u$ be a solution to $\square u=F$ on a strip $[a, b]\times\R^3$.
\begin{eqnarray}\label{waveop-esti}
   \lefteqn{{\|u\|}_{L_t^{\frac{2}{\gamma}} L_x^{\frac{2}{1-\gamma}}([a,\,b]\times\R^3)}
   +{\|u\|}_{C_t\dot{H}_x^\gamma([a,\,b]\times\R^3)}
   +{\|\partial_t u\|}_{C_t\dot{H}_x^{\gamma-1}([a,\,b]\times\R^3)}}
   \nonumber\\ & \le & C^\ast_\gamma\,\bigg({\|u(a)\|}_{\dot{H}_x^\gamma(\R^3)}
   +{\|\partial_t u(a)\|}_{\dot{H}_x^{\gamma-1}(\R^3)}
   +{\|F\|}_{L_t^{\frac{2}{1+\gamma}} L_x^{\frac{2}{2-\gamma}}([a,\,b]\times\R^3)}\bigg).
\end{eqnarray}
The constant $C^\ast_\gamma$ is independent of $a$ and $b$.

Next an estimate for ${\cal W}$ is derived. It might not be optimal,
but it will turn out to be sufficient in the sequel.

\begin{lemma} Let ${\cal W}$ be defined by (\ref{calW-def}). If $T>0$ and $S_T=[0, T]\times\R^3$, then
\begin{equation}\label{calW-esti1}
   {\|{\cal W}h\|}_{L^\infty_t\dot{H}^{1-\eps}_x(S_T)}
   \le C\eps^{-1}T^\eps\,{\|h\|}_{L^\infty_t L^2_x(S_T)},
   \quad\eps\in ]0, 1].
\end{equation}
In particular,
\begin{equation}\label{calW-esti2}
   {\|{\cal W}h\|}_{L^\infty_t L^{\frac{6}{1+2\eps}}_x(S_T)}
   \le C_2(\eps, T)\,{\|h\|}_{L^\infty_t L^2_x(S_T)},
\end{equation}
where $C_2(\eps, T)=C_1(\eps)C\eps^{-1}T^\eps$ is increasing in $T$.
\end{lemma}
{\bf Proof\,:} The Fourier transform of ${\cal W}h$ is
\begin{eqnarray*}
   \widehat{({\cal W}h)}(t, \xi) & = & \int_{\R^3} e^{-i\,\xi\cdot x}\,({\cal W}h)(t, x)\,dx
   =\int_{|y|\le t}\frac{dy}{4\pi |y|^2}\,\int_{\R^3} dx\,e^{-i\,\xi\cdot x}\,h(t-|y|, x+y)
   \\ & = & \int_{|y|\le t}\frac{dy}{4\pi |y|^2}\,e^{i\,\xi\cdot y}\,\hat{h}(t-|y|, \xi)
   =\int_0^t\frac{ds}{4\pi}\,\hat{h}(t-s, \xi)\,\int_{|\omega|=1} dS(\omega)\,e^{is\,\xi\cdot\omega}
   \\ & = & \int_0^t\frac{\sin(s|\xi|)}{s|\xi|}\,\hat{h}(t-s, \xi)\,ds.
\end{eqnarray*}
Now use $|\sin(s|\xi|)|\le\min\{1, s|\xi|\}$ to obtain for $\eps\in ]0, 1[$
\begin{eqnarray*}
   |\widehat{({\cal W}h)}(t, \xi)|
   & \le & \int_0^t {\bf 1}_{\{1\le s|\xi|\}}\frac{1}{s|\xi|}\,|\hat{h}(t-s, \xi)|\,ds
   +\int_0^t {\bf 1}_{\{1<\frac{1}{s|\xi|}\}}\,|\hat{h}(t-s, \xi)|\,ds
   \\ & \le & \int_0^t {\bf 1}_{\{1\le s|\xi|\}}\frac{1}{s|\xi|}\,(s|\xi|)^\eps\,|\hat{h}(t-s, \xi)|\,ds
   +\int_0^t {\bf 1}_{\{1<\frac{1}{s|\xi|}\}}\,\bigg(\frac{1}{s|\xi|}\bigg)^{1-\eps}\,|\hat{h}(t-s, \xi)|\,ds
   \\ & \le & \frac{2}{|\xi|^{1-\eps}}\int_0^t\frac{ds}{s^{1-\eps}}\,|\hat{h}(t-s, \xi)|.
\end{eqnarray*}
It follows that
\begin{eqnarray*}
   {\|({\cal W}h)(t)\|}_{\dot{H}^{1-\eps}_x(\R^3)}
   & = & \frac{1}{(2\pi)^{3/2}}\,{\||\xi|^{1-\eps}\,\widehat{({\cal W}h)}(t)\|}_{L^2_\xi(\R^3)}
   \\ & \le & 2\int_0^t\frac{ds}{s^{1-\eps}}\,{\|h(t-s)\|}_{L^2_x(\R^3)},
\end{eqnarray*}
and this yields (\ref{calW-esti1}). By the homogeneous Sobolev embedding in $\R^3$,
\[ {\|{\cal W}h\|}_{L^\infty_t L^{\frac{6}{1+2\eps}}_x(S_T)}
   \le C_1(\eps) {\|{\cal W}h\|}_{L^\infty_t\dot{H}^{1-\eps}_x(S_T)}
   \le C_1(\eps)C\eps^{-1}T^\eps\,{\|h\|}_{L^\infty_t L^2_x(S_T)}, \]
which is (\ref{calW-esti2}). {\hfill$\Box$}\bigskip

\begin{cor}\label{flat-cor} For $\eps\in ]0, 1]$,
\[ {\Big\||E_\flat|+|B_\flat|\Big\|}_{L^\infty_t L^{\frac{6}{1+2\eps}}_x(S_T)}
   \le C_3(\eps, T)\,{\|\sigma_{-1}\|}_{L^\infty_t L^2_x(S_T)}, \]
where $C_3$ is increasing in $T$. In particular,
\begin{equation}\label{efus-esti3}
   {\|\,\square^{-1}((|E_\flat|+|B_\flat|)\sigma_{-1})\|}_{L^\frac{3}{1-\eps}_t L^\frac{6}{1+2\eps}_x(S_T)}
   \le C_4(\eps, T)\,{\|\sigma_{-1}\|}_{L^\infty_t L^2_x(S_T)}^2
\end{equation}
for a constant $C_4$ that is increasing in $T$.
\end{cor}
{\bf Proof\,:} The first bound is an immediate consequence
of (\ref{EB-T-esti}) and (\ref{calW-esti2}). To prove (\ref{efus-esti3}),
note that
\begin{eqnarray*}
   {\|(|E_\flat|+|B_\flat|)\sigma_{-1}\|}_{L^\infty_t L^\frac{3}{2+\eps}_x(S_T)}
   & \le & {\Big\||E_\flat|+|B_\flat|\Big\|}_{L^\infty_t L^\frac{6}{1+2\eps}_x(\R^3)}
   {\|\sigma_{-1}\|}_{L^\infty_t L^2_x(\R^3)}
   \\ & \le & C_3(\eps, T)\,{\|\sigma_{-1}\|}_{L^\infty_t L^2_x(S_T)}^2.
\end{eqnarray*}
Thus using (\ref{waveop-esti}) for $\gamma_\eps=\frac{2}{3}(1-\eps)$,
where $\frac{2}{2-\gamma_\eps}=\frac{3}{2+\eps}$, $\frac{2}{1+\gamma_\eps}=\frac{6}{5-2\eps}$,
$\frac{2}{1-\gamma_\eps}=\frac{6}{1+2\eps}$, and $\frac{2}{\gamma_\eps}=\frac{3}{1-\eps}$,
\begin{eqnarray*}
   {\|\,\square^{-1}((|E_\flat|+|B_\flat|)\sigma_{-1})\|}_{L^\frac{3}{1-\eps}_t L^\frac{6}{1+2\eps}_x(S_T)}
   & \le & C^\ast_{\gamma_\eps}\,
   {\|(|E_\flat|+|B_\flat|)\sigma_{-1}\|}_{L^{\frac{6}{5-2\eps}}_t L^\frac{3}{2+\eps}_x(S_T)}
   \\ & \le & C^\ast_{\gamma_\eps}\,
   T^{\frac{5-2\eps}{6}}\,C_3(\eps, T)\,{\|\sigma_{-1}\|}_{L^\infty_t L^2_x(S_T)}^2,
\end{eqnarray*}
as was to be shown. {\hfill$\Box$}\bigskip

\begin{lemma}\label{sharp-lem} For $\eps\in ]0, 1]$,
\[ {\Big\||E_\sharp|+|B_\sharp|\Big\|}_{L^\frac{3}{1-\eps}_t L^\frac{6}{1+2\eps}_x(S_T)}
   \le C_5(\eps, T, {\rm data}, {\|\sigma_{-1}\|}_{L^\infty_t L^2_x(S_T)}), \]
where $C_5$ is increasing in both the $T$-argument and the $\|\cdot\|$-argument.
\end{lemma}
{\bf Proof\,:} The argument is similar to \cite[Thm.~4.8, p.~108]{sogge}.
Fix an interval $[a, b]\subset [0, T]$. According to (\ref{EB-S-esti}),
(\ref{E-form}) and (\ref{B-form}), and (\ref{efus-esti3}),
\begin{eqnarray*}
   \lefteqn{{\Big\||E_\sharp|+|B_\sharp|\Big\|}_{L^\frac{3}{1-\eps}_t L^\frac{6}{1+2\eps}_x([a,\,b]\times\R^3)}}
   \\ & \le & {\|\,\square^{-1}((|E|+|B|)
   \sigma_{-1})\|}_{L^\frac{3}{1-\eps}_t L^\frac{6}{1+2\eps}_x([a,\,b]\times\R^3)}
   \\ & \le & {\|\,\square^{-1}((|E_D|+|E_{DT}|
   +|B_D|+|B_{DT}|)\sigma_{-1})\|}_{L^\frac{3}{1-\eps}_t L^\frac{6}{1+2\eps}_x([a,\,b]\times\R^3)}
   \\ & & +\,{\|\,\square^{-1}((|E_\flat|+|B_\flat|)
   \sigma_{-1})\|}_{L^\frac{3}{1-\eps}_t L^\frac{6}{1+2\eps}_x(S_T)}
   +{\|\,\square^{-1}((|E_\sharp|+|B_\sharp|)
   \sigma_{-1})\|}_{L^\frac{3}{1-\eps}_t L^\frac{6}{1+2\eps}_x([a,\,b]\times\R^3)}
   \\ & \le & ({\rm data})+C_4(\eps, T)\,{\|\sigma_{-1}\|}_{L^\infty_t L^2_x(S_T)}^2
   +\,{\|\,\square^{-1}((|E_\sharp|+|B_\sharp|)
   \sigma_{-1})\|}_{L^\frac{3}{1-\eps}_t L^\frac{6}{1+2\eps}_x([a,\,b]\times\R^3)}.
\end{eqnarray*}
Let $F=(|E_\sharp|+|B_\sharp|)\sigma_{-1}$ and $u=\square^{-1}F$.
Then (\ref{waveop-esti}) for $\gamma_\eps=\frac{2}{3}(1-\eps)$ yields
\begin{eqnarray*}
   \lefteqn{{\|u\|}_{L^\frac{3}{1-\eps}_t L^\frac{6}{1+2\eps}_x([a,\,b]\times\R^3)}}
   \\ & \le & C^\ast_{\gamma_\eps}\,\bigg({\|u(a)\|}_{\dot{H}_x^{\gamma_\eps}(\R^3)}
   +{\|\partial_t u(a)\|}_{\dot{H}_x^{\gamma_\eps-1}(\R^3)}
   +{\|F\|}_{L^{\frac{6}{5-2\eps}}_t L^\frac{3}{2+\eps}_x([a,\,b]\times\R^3)}\bigg)
   \\ & \le & C^\ast_{\gamma_\eps}\,\bigg({\|u(a)\|}_{\dot{H}_x^{\gamma_\eps}(\R^3)}
   +{\|\partial_t u(a)\|}_{\dot{H}_x^{\gamma_\eps-1}(\R^3)}
   \\ & & \hspace{3em} +\,{\Big\||E_\sharp|
   +|B_\sharp|\Big\|}_{L^{\frac{3}{1-\eps}}_t L^\frac{6}{1+2\eps}_x([a,\,b]\times\R^3)}
   {\|\sigma_{-1}\|}_{L^2_t L^2_x([a,\,b]\times\R^3)}\bigg).
\end{eqnarray*}
As a consequence,
\begin{eqnarray}\label{detojeh}
   \lefteqn{{\Big\||E_\sharp|
   +|B_\sharp|\Big\|}_{L^\frac{3}{1-\eps}_t L^\frac{6}{1+2\eps}_x([a,\,b]\times\R^3)}}
   \nonumber\\ & \le & ({\rm data})+C_4(\eps, T)\,{\|\sigma_{-1}\|}_{L^\infty_t L^2_x(S_{T})}^2
   +\,C^\ast_{\gamma_\eps}\,\bigg({\|u(a)\|}_{\dot{H}_x^{\gamma_\eps}(\R^3)}
   +{\|\partial_t u(a)\|}_{\dot{H}_x^{\gamma_\eps-1}(\R^3)}
   \nonumber\\ & & \hspace{17em} +\,{\Big\||E_\sharp|
   +|B_\sharp|\Big\|}_{L^{\frac{3}{1-\eps}}_t L^\frac{6}{1+2\eps}_x([a,\,b]\times\R^3)}
   {\|\sigma_{-1}\|}_{L^2_t L^2_x([a,\,b]\times\R^3)}\bigg).
   \nonumber\\ & &
\end{eqnarray}
Without loss of generality suppose that $C^\ast_{\gamma_\eps}
{\|\sigma_{-1}\|}_{L^2_t L^2_x(S_T)}\ge 1$, since otherwise one can just take
$[a, b]=[0, T]$. Fix a finite partition $0=T_0<T_1<T_2<\ldots<T_{N-1}<T_N=T$
of $[0, T]$ such that
\begin{equation}\label{sig-1-teil}
   {\|\sigma_{-1}\|}_{L^2_t L^2_x([T_j,\,T_{j+1}])}=\frac{1}{2C^\ast_{\gamma_\eps}}
   \,\,\,(j=0, \ldots, N-2),\quad {\|\sigma_{-1}\|}_{L^2_t L^2_x([T_{N-1},\,T_N])}
   \le\frac{1}{2C^\ast_{\gamma_\eps}}.
\end{equation}
Note that
\[ (N-1)\frac{1}{4(C^\ast_{\gamma_\eps})^2}
   \le\sum_{j=0}^{N-2}{\|\sigma_{-1}\|}_{L^2_t L^2_x([T_{N-1},\,T_N])}^2
   \le {\|\sigma_{-1}\|}_{L^2_t L^2_x(S_T)}^2 \]
yields the upper bound
\begin{equation}\label{sgNbd}
   N\le 1+4(C^\ast_{\gamma_\eps})^2{\|\sigma_{-1}\|}_{L^2_t L^2_x(S_T)}^2.
\end{equation}
By (\ref{detojeh}) and (\ref{sig-1-teil}),
\begin{eqnarray}\label{flesup}
   \lefteqn{{\Big\||E_\sharp|
   +|B_\sharp|\Big\|}_{L^\frac{3}{1-\eps}_t L^\frac{6}{1+2\eps}_x([T_j,\,T_{j+1}]\times\R^3)}}
   \nonumber\\ & \le & 2({\rm data})_j+2C_4(\eps, T)\,{\|\sigma_{-1}\|}_{L^\infty_t L^2_x(S_{T})}^2
   +\,2C^\ast_{\gamma_\eps}\,\bigg({\|u(T_j)\|}_{\dot{H}_x^{\gamma_\eps}(\R^3)}
   +{\|\partial_t u(T_j)\|}_{\dot{H}_x^{\gamma_\eps-1}(\R^3)}\bigg). \qquad
\end{eqnarray}
Thus in particular
\begin{eqnarray*}
   {\|F\|}_{L^{\frac{6}{5-2\eps}}_t L^\frac{3}{2+\eps}_x([T_j,\,T_{j+1}]\times\R^3)}
   &\le & \frac{1}{2C^\ast_{\gamma_\eps}}\,{\Big\||E_\sharp|+|B_\sharp|\Big\|}_{L^\frac{3}{1-\eps}_t
   L^\frac{6}{1+2\eps}_x([T_j,\,T_{j+1}]\times\R^3)}
   \\ & \le & \frac{1}{C^\ast_{\gamma_\eps}}\,({\rm data})_j
   +\frac{1}{C^\ast_{\gamma_\eps}}\,C_4(\eps, T)\,{\|\sigma_{-1}\|}_{L^\infty_t L^2_x(S_{T})}^2
   \\ & & +\,{\|u(T_j)\|}_{\dot{H}_x^{\gamma_\eps}(\R^3)}
   +{\|\partial_t u(T_j)\|}_{\dot{H}_x^{\gamma_\eps-1}(\R^3)},
\end{eqnarray*}
so that by (\ref{waveop-esti}) for the interval $[T_j, T_{j+1}]$,
\begin{eqnarray*}
   \lefteqn{{\|u(T_{j+1})\|}_{\dot{H}_x^{\gamma_\eps}(\R^3)}
   +{\|\partial_t u(T_{j+1})\|}_{\dot{H}_x^{\gamma_\eps-1}(\R^3)}}
   \\ & \le & C^\ast_{\gamma_\eps}\,\bigg({\|u(T_j)\|}_{\dot{H}_x^{\gamma_\eps}(\R^3)}
   +{\|\partial_t u(T_j)\|}_{\dot{H}_x^{\gamma_\eps-1}(\R^3)}
   +{\|F\|}_{L^{\frac{6}{5-2\eps}}_t L^\frac{3}{2+\eps}_x([T_j,\,T_{j+1}]\times\R^3)}\bigg)
   \\ & \le & ({\rm data})_j+C_4(\eps, T)\,{\|\sigma_{-1}\|}_{L^\infty_t L^2_x(S_{T})}^2
   +2C^\ast_{\gamma_\eps}\,\bigg({\|u(T_j)\|}_{\dot{H}_x^{\gamma_\eps}(\R^3)}
   +{\|\partial_t u(T_j)\|}_{\dot{H}_x^{\gamma_\eps-1}(\R^3)}\bigg). \qquad
\end{eqnarray*}
Iteration of this estimate and noting that $u(0)=\partial_t u(0)=0$
leads to the bound
\[ {\|u(T_j)\|}_{\dot{H}_x^{\gamma_\eps}(\R^3)}
   +{\|\partial_t u(T_j)\|}_{\dot{H}_x^{\gamma_\eps-1}(\R^3)}
   \le\sum_{i=0}^{j-1} (2C^\ast_{\gamma_\eps})^{j-1-i}\Big(
   ({\rm data})_i+C_4(\eps, T)\,{\|\sigma_{-1}\|}_{L^\infty_t L^2_x(S_{T})}^2\Big). \]
Hence by (\ref{flesup}),
\begin{eqnarray*}
   \lefteqn{{\Big\||E_\sharp|
   +|B_\sharp|\Big\|}_{L^\frac{3}{1-\eps}_t L^\frac{6}{1+2\eps}_x([T_j,\,T_{j+1}]\times\R^3)}}
   \\ & \le & 2\sum_{i=0}^j (2C^\ast_{\gamma_\eps})^{j-i} ({\rm data})_i
   +2C_4(\eps, T)\,\bigg(\sum_{i=0}^j (2C^\ast_{\gamma_\eps})^i\bigg)
   \,{\|\sigma_{-1}\|}_{L^\infty_t L^2_x(S_T)}^2
\end{eqnarray*}
for $j=0, \ldots, N-1$. Therefore the estimate
\begin{eqnarray*}
   \lefteqn{{\Big\||E_\sharp|+|B_\sharp|\Big\|}_{L^\frac{3}{1-\eps}_t L^\frac{6}{1+2\eps}_x(S_T)}}
   \\ & \le & \sum_{j=0}^{N-1}\,{\Big\||E_\sharp|+|B_\sharp|\Big\|}_{L^\frac{3}{1-\eps}_t L^\frac{6}{1+2\eps}_x
   ([T_j,\,T_{j+1}]\times\R^3)}
   \\ & \le & 2\sum_{j=0}^{N-1}\sum_{i=0}^j (2C^\ast_{\gamma_\eps})^{j-i} ({\rm data})_i
   +2C_4(\eps, T)\,\bigg(\sum_{j=0}^{N-1}\sum_{i=0}^j (2C^\ast_{\gamma_\eps})^i\bigg)
   \,{\|\sigma_{-1}\|}_{L^\infty_t L^2_x(S_T)}^2
\end{eqnarray*}
is obtained from (\ref{sgNbd}). Recalling (\ref{sgNbd}), the claim is obtained.
{\hfill$\Box$}\bigskip

\begin{cor}\label{E+B} For $\eps\in ]0, 1]$,
\[ {\Big\||E|+|B|\Big\|}_{L^\frac{3}{1-\eps}_t L^\frac{6}{1+2\eps}_x(S_T)}
   \le C_6(\eps, T, {\rm data}, {\|\sigma_{-1}\|}_{L^\infty_t L^2_x(S_T)}), \]
where $C_6$ is increasing in both the $T$-argument and the $\|\cdot\|$-argument.
\end{cor}
{\bf Proof\,:} From (\ref{E-form}) and (\ref{B-form}),
\[ E=E_D+E_{DT}+E_\flat+E_\sharp\quad\mbox{and}\quad
   B=B_D+B_{DT}+B_\flat+B_\sharp, \]
where $E_D$, $E_{DT}$, $B_D$, and $B_{DT}$ are data terms.
Since Corollary \ref{flat-cor} implies that
\[ {\Big\||E_\flat|+|B_\flat|\Big\|}_{L^\frac{3}{1-\eps}_t L^{\frac{6}{1+2\eps}}_x(S_T)}
   \le C_3(\eps, T)\,T^{\frac{1-\eps}{3}}\,{\|\sigma_{-1}\|}_{L^\infty_t L^2_x(S_T)}, \]
it remains to apply Lemma \ref{sharp-lem}. {\hfill$\Box$}\bigskip

\begin{cor}\label{mom-cor} If $\eps\in ]0, \frac{1}{10}]$, then
\[ m_{\frac{3(1-2\eps)}{1+2\eps}}(t)\le
   C_7(0, \eps, t, {\rm data}, {\|\sigma_{-1}\|}_{L^\infty_t L^2_x(S_t)}), \]
where $C_7$ is increasing in both the $t$-argument and the $\|\cdot\|$-argument.
\end{cor}
{\bf Proof\,:} By Lemma \ref{moment-ode} for $k=\frac{3(1-2\eps)}{1+2\eps}\ge 2$,
\[ m_{\frac{3(1-2\eps)}{1+2\eps}}(t)\le m_{\frac{3(1-2\eps)}{1+2\eps}}(0)
   +C(0, \eps)\,\int_0^t {\|E(s)\|}_{L^{\frac{6}{1+2\eps}}_x(\R^3)}
   \,m_{\frac{3(1-2\eps)}{1+2\eps}}(s)^{\frac{5-2\eps}{6}}\,ds. \]
By a standard differential inequality comparison theorem, this yields
\[ m_{\frac{3(1-2\eps)}{1+2\eps}}(t)^{\frac{1+2\eps}{6}}
   \le m_{\frac{3(1-2\eps)}{1+2\eps}}(0)^{\frac{1+2\eps}{6}}
   +C(0, \eps)\bigg(\frac{1+2\eps}{6}\bigg)
   \int_0^t {\|E(s)\|}_{L^{\frac{6}{1+2\eps}}_x(\R^3)}\,ds, \]
so that by Corollary \ref{E+B},
\begin{eqnarray*}
   m_{\frac{3(1-2\eps)}{1+2\eps}}(t)
   & \le & C(0, \eps, {\rm data})\bigg(
   1+\bigg[\int_0^t {\|E(s)\|}_{L^{\frac{6}{1+2\eps}}_x(\R^3)}
   \,ds\bigg]^{\frac{6}{1+2\eps}}\bigg)
   \\ & \le & C(0, \eps, {\rm data})\bigg(
   1+t^{\frac{2(2+\eps)}{1+2\eps}}
   \,{\|E\|}_{L^\frac{3}{1-\eps}_t L^{\frac{6}{1+2\eps}}_x(S_t)}^{\frac{6}{1+2\eps}}\bigg)
   \\ & \le & C(0, \eps, {\rm data})\bigg(
   1+t^{\frac{2(2+\eps)}{1+2\eps}}
   \,C_6(\eps, t, {\rm data}, {\|\sigma_{-1}\|}_{L^\infty_t L^2_x(S_t)})^{\frac{6}{1+2\eps}}\bigg).
\end{eqnarray*}
Hence $C_7$ can be defined appropriately. {\hfill$\Box$}\bigskip

\begin{cor}\label{sig-1L4} If $\eps\in ]0, \frac{1}{10}]$, then
\[ {\|\sigma_{-1}(t)\|}_{L^{\frac{4(1-\eps)}{1+2\eps}}_x(\R^3)}
   \le C_8(0, \eps, t, {\rm data}, {\|\sigma_{-1}\|}_{L^\infty_t L^2_x(S_t)})
   \,P_\infty(t), \]
where $C_8$ is increasing in both the $t$-argument and the $\|\cdot\|$-argument.
\end{cor}
{\bf Proof\,:} In Lemma \ref{sig-1einP} take $q=\frac{4(1-\eps)}{1+2\eps}$
and $\alpha=\frac{3(1-2\eps)}{8(1-\eps)}$. Invoking Corollary \ref{mom-cor},
it follows that
\begin{eqnarray*}
   {\|\sigma_{-1}(t)\|}_{L^{\frac{4(1-\eps)}{1+2\eps}}_x(\R^3)}
   & \le & C(0, \eps)\,P_\infty(t)
   \,m_{\frac{3(1-2\eps)}{1+2\eps}}(t)^{\frac{1+2\eps}{4(1-\eps)}}
   \\ & \le & C(0, \eps)\,P_\infty(t)
   \,C_7(0, \eps, t, {\rm data}, {\|\sigma_{-1}\|}_{L^\infty_t L^2_x(S_t)})
   ^{\frac{1+2\eps}{4(1-\eps)}}.
\end{eqnarray*}
Thus it remains to choose $C_8$ in a suitable manner.
{\hfill$\Box$}\bigskip


\setcounter{equation}{0}

\section{Proof of Theorem \ref{sig-1thm}}
\label{mainprf-sect} 

Fix any characteristic $(X_0, P_0)$ in the support of $f$, i.e.,
\begin{equation}\label{charinsupp}
   f(s, X_0(s), P_0(s))=f^{(0)}(X_0(0), P_0(0))\neq 0,\quad s\in [0, T_{{\rm max}}[.
\end{equation}
The relation
\begin{eqnarray*}
   \frac{d}{ds}\,\sqrt{1+P_0(s)^2}
   & = & V_0(s)\cdot\dot{P}_0(s)
   =V_0(s)\cdot\Big(E(s, X_0(s))+V_0(s)\wedge B(s, X_0(s))\Big)
   \\ & = & V_0(s)\cdot E(s, X_0(s))
\end{eqnarray*}
in conjunction with (\ref{E-form}) yields
\begin{eqnarray}\label{P0bd}
   |P_0(t)| & \le & \sqrt{1+P_0(t)^2}
   =\sqrt{1+P_0(0)^2}+\int_0^t V_0(s)\cdot E(s, X_0(s))\,ds
   \nonumber \\ & = & \sqrt{1+P_0(0)^2}+\int_0^t V_0(s)\cdot (E_D+E_{DT})(s, X_0(s))\,ds
   \nonumber \\ & & +\int_0^t V_0(s)\cdot E_\flat(s, X_0(s))\,ds
   +\int_0^t V_0(s)\cdot E_\sharp(s, X_0(s))\,ds.
\end{eqnarray}
By the definition of $E_\flat$,
\begin{eqnarray*}
   I_\flat(t) & = & \int_0^t V_0(s)\cdot E_\flat(s, X_0(s))\,ds
   \\ & = & -\int_0^t ds\,\int_{|y|\le s}\frac{dy}{|y|^2}
   \int_{\R^3}\frac{dp}{1+p^2}\,\frac{1}{(1+v\cdot\omega)^2}
   \,V_0(s)\cdot (v+\omega)\,f(s-|y|, X_0(s)+y, p)
   \\ & = & -\int_0^t d\tau\,\int_\tau^t ds\,\int_{|\omega|=1} dS(\omega)
   \\ & & \hspace{7em}\times\int_{\R^3}\frac{dp}{1+p^2}\,\frac{1}{(1+v\cdot\omega)^2}
   \,V_0(s)\cdot (v+\omega)\,f(\tau, X_0(s)+(s-\tau)\omega, p).
\end{eqnarray*}
Next write $V_0(s)\cdot (v+\omega)=(V_0(s)+\omega)\cdot (v+\omega)
-(1+v\cdot\omega)$ and split the integral accordingly
as $I_\flat(t)=I_{\flat, 1}(t)+I_{\flat, 2}(t)$. Firstly,
\begin{eqnarray*}
   |I_{\flat, 2}(t)| & = & \int_0^t d\tau\,\int_\tau^t ds\,\int_{|\omega|=1} dS(\omega)
   \int_{\R^3}\frac{dp}{1+p^2}\,\frac{1}{1+v\cdot\omega}
   \,f(\tau, X_0(s)+(s-\tau)\omega, p)
   \\ & \le & {\cal L}(0)\int_0^t d\tau\,\int_\tau^t ds\,\int_{|\omega|=1} dS(\omega)
   \int_{|p|\le P_\infty(\tau)}\frac{dp}{1+p^2}\,\frac{1}{1+v\cdot\omega}
   \\ & \le & C(0)\int_0^t (t-\tau)\,\ln P_\infty(\tau)\,P_\infty(\tau)\,d\tau
   \le C(0)\,t\int_0^t\ln P_\infty(\tau)\,P_\infty(\tau)\,d\tau, 
\end{eqnarray*}
where we used Lemma \ref{efu-lem1}(a). Concerning $I_{\flat, 1}(t)$,
\begin{eqnarray*}
   |I_{\flat, 1}(t)| & = & \bigg|\int_0^t d\tau\,\int_\tau^t ds\,\int_{|\omega|=1} dS(\omega)
   \\ & & \hspace{4em}\times\int_{\R^3}\frac{dp}{1+p^2}\,
   \frac{1}{(1+v\cdot\omega)^2}\,(V_0(s)+\omega)\cdot (v+\omega)
   \,f(\tau, X_0(s)+(s-\tau)\omega, p)\bigg|
   \\ & \le & \int_0^t d\tau\,\int_\tau^t ds\,\int_0^{2\pi} d\varphi\int_0^\pi d\theta\,\sin\theta
   \\ & & \hspace{4em}\times\int_{\R^3}\frac{dp}{1+p^2}\,
   \frac{(1+V_0(s)\cdot\omega)^{1/2}}{(1+v\cdot\omega)^{3/2}}
   \,f(\tau, X_0(s)+(s-\tau)\omega, p)
\end{eqnarray*}
for $\omega=(\cos\varphi\sin\theta, \sin\varphi\sin\theta, \cos\theta)$;
recall (\ref{vecnorm1}). If $tP_\infty(t)\ge 1$, then the $\int_0^t d\tau\int_\tau^t ds$
is split to find
\[ {\bf 1}_{\{tP_\infty(t)\ge 1\}}\,|I_{\flat, 1}(t)|
   \le I_{\flat, 11}(t)+I_{\flat, 12}(t)+I_{\flat, 13}(t), \]
where
\begin{eqnarray*}
   I_{\flat, 11}(t) & = & \int_0^{t-P_\infty(t)^{-1}} d\tau\,\int_\tau^{\tau+P_\infty(t)^{-1}} ds\,(\ldots),
   \\ I_{\flat, 12}(t) & = & \int_0^{t-P_\infty(t)^{-1}} d\tau\,\int_{\tau+P_\infty(t)^{-1}}^t ds\,(\ldots),
   \\ I_{\flat, 13}(t) & = & \int_{t-P_\infty(t)^{-1}}^t d\tau\,\int_\tau^t ds\,(\ldots).
\end{eqnarray*}
To begin with, by Lemma \ref{efu-lem1}(b) for $\theta=1$ and $\kappa=\frac{3}{2}$,
\begin{eqnarray*}
   I_{\flat, 11}(t) & = & \int_0^{t-P_\infty(t)^{-1}} d\tau\,\int_\tau^{\tau+P_\infty(t)^{-1}} ds
   \,\int_0^{2\pi} d\varphi\int_0^\pi d\theta\,\sin\theta
   \\ & & \hspace{4em}\times\int_{\R^3}\frac{dp}{1+p^2}\,
   \frac{(1+V_0(s)\cdot\omega)^{1/2}}{(1+v\cdot\omega)^{3/2}}
   \,f(\tau, X_0(s)+(s-\tau)\omega, p)
   \\ & \le & C{\cal L}(0)\,P_\infty(t)^{-1}\int_0^{t-P_\infty(t)^{-1}} d\tau
   \,\max_{|\omega|=1}\int_{|p|\le P_\infty(\tau)}\frac{dp}{1+p^2}\,
   \frac{1}{(1+v\cdot\omega)^{3/2}}
   \\ & \le & C(0)\,P_\infty(t)^{-1}\int_0^t P_\infty(\tau)^2\,d\tau
   \le C(0)\,\int_0^t P_\infty(\tau)\,d\tau.
\end{eqnarray*}
Note that in the last step it was used that $P_\infty(\tau)\le P_\infty(t)$,
since $P_\infty$ is increasing by definition. Similarly,
\begin{eqnarray*}
   I_{\flat, 13}(t) & = & \int_{t-P_\infty(t)^{-1}}^t d\tau\,\int_\tau^t ds
   \,\int_0^{2\pi} d\varphi\int_0^\pi d\theta\,\sin\theta
   \\ & & \hspace{4em}\times\int_{\R^3}\frac{dp}{1+p^2}\,
   \frac{(1+V_0(s)\cdot\omega)^{1/2}}{(1+v\cdot\omega)^{3/2}}
   \,f(\tau, X_0(s)+(s-\tau)\omega, p)
   \\ & \le & {\cal L}(0)\int_{t-P_\infty(t)^{-1}}^t d\tau\,(t-\tau)
   \,\max_{|\omega|=1}\int_{|p|\le P_\infty(\tau)}\frac{dp}{1+p^2}\,
   \frac{1}{(1+v\cdot\omega)^{3/2}}
   \\ & \le & C{\cal L}(0)\int_{t-P_\infty(t)^{-1}}^t d\tau\,(t-\tau)
   \,P_\infty(\tau)^2
   \\ & \le & C(0)\,P_\infty(t)^{-1}\int_0^t P_\infty(\tau)^2\,d\tau
   \le C(0)\int_0^t P_\infty(\tau)\,d\tau.
\end{eqnarray*}
It remains to deal with
\begin{eqnarray*}
   I_{\flat, 12}(t) & = & \int_0^{t-P_\infty(t)^{-1}} d\tau\,\int_{\tau+P_\infty(t)^{-1}}^t ds
   \,\int_0^{2\pi} d\varphi\int_0^\pi d\theta\,\sin\theta
   \\ & & \hspace{4em}\times\int_{\R^3}\frac{dp}{1+p^2}\,
   \frac{(1+V_0(s)\cdot\omega)^{1/2}}{(1+v\cdot\omega)^{3/2}}
   \,f(\tau, X_0(s)+(s-\tau)\omega, p).
\end{eqnarray*}
Let $M_\tau=[\tau+P_\infty(t)^{-1}, t]\times [0, 2\pi]\times [0, \pi]$
and consider the mapping
\begin{equation}\label{palla-trans}
   \Phi_\tau:\quad M_\tau\ni (s, \varphi, \theta)\mapsto y=X_0(s)+(s-\tau)\omega\in\R^3.
\end{equation}
According to \cite[Lemma 2.1]{pallard}, $\Phi_\tau$ is a diffeomorphism and such that
\[ dy=(1+V_0(s)\cdot\omega)(s-\tau)^2\sin\theta\,ds\,d\varphi\,d\theta. \]
Writing the inverse mapping as $s=s(y)$ and $\omega=\omega(y)$,
this yields using Lemma \ref{efu-lem1}(b) for $\theta=\kappa=3$,
\begin{eqnarray*}
   \lefteqn{I_{\flat, 12}(t)}
   \\ & = & \int_0^{t-P_\infty(t)^{-1}} d\tau\,\int_{\Phi_\tau(M_\tau)} dy\,\frac{1}{(s-\tau)^2}
   \,\frac{1}{(1+V_0(s)\cdot\omega)^{1/2}}
   \\ & & \hspace{11em}\times\int_{|p|\le P_\infty(\tau)}\frac{dp}{1+p^2}
   \,\frac{1}{(1+v\cdot\omega)^{3/2}}\,f(\tau, y, p)
   \\ & \le & \int_0^{t-P_\infty(t)^{-1}} d\tau\,
   \bigg(\int_{\Phi_\tau(M_\tau)} dy\,
   \frac{1}{(s-\tau)^4}\,\frac{1}{(1+V_0(s)\cdot\omega)}
   \int_{|p|\le P_\infty(\tau)}
   \,\frac{dp}{(1+p^2)^3}\,\frac{1}{(1+v\cdot\omega)^3}\bigg)^{1/2}
   \\ & & \hspace{7em} \times\bigg(\int_{\Phi_\tau(M_\tau)}\int_{|p|\le P_\infty(\tau)} dy\,dp\,(1+p^2)
   \,f(\tau, y, p)^2\bigg)^{1/2}
   \\ & \le & C{\cal L}(0)^{1/2}\int_0^{t-P_\infty(t)^{-1}} d\tau\,
   \bigg(\int_{\Phi_\tau(M_\tau)} dy\,
   \frac{1}{(s-\tau)^4}\,\frac{1}{(1+V_0(s)\cdot\omega)}\bigg)^{1/2}
   \,P_\infty(\tau)^{1/2}\,m_2(\tau)^{1/2}
   \\ & \le & C(0)\int_0^{t-P_\infty(t)^{-1}} d\tau\,
   \bigg(\int_{\tau+P_\infty(t)^{-1}}^t ds
   \,\int_0^{2\pi} d\varphi\int_0^\pi d\theta\,\sin\theta\,
   \frac{1}{(s-\tau)^2}\bigg)^{1/2}
   \,P_\infty(\tau)^{1/2}\,m_2(\tau)^{1/2}
   \\ & \le & C(0)P_\infty(t)^{1/2}\int_0^t
   \,P_\infty(\tau)^{1/2}\,m_2(\tau)^{1/2}\,d\tau.
\end{eqnarray*}
To summarize, it has been shown that
\begin{eqnarray}\label{Iflat-bd}
   {\bf 1}_{\{tP_\infty(t)\ge 1\}}\,|I_{\flat}(t)|
   & \le & {\bf 1}_{\{tP_\infty(t)\ge 1\}}\,|I_{\flat, 11}(t)|
   +{\bf 1}_{\{tP_\infty(t)\ge 1\}}\,|I_{\flat, 12}(t)|
   +{\bf 1}_{\{tP_\infty(t)\ge 1\}}\,|I_{\flat, 13}(t)|
   \nonumber \\ & & +\,{\bf 1}_{\{tP_\infty(t)\ge 1\}}\,|I_{\flat, 2}(t)|
   \nonumber \\ & \le & C(0)\int_0^t P_\infty(\tau)\,d\tau
   +C(0)P_\infty(t)^{1/2}\int_0^t
   \,P_\infty(\tau)^{1/2}\,m_2(\tau)^{1/2}\,d\tau
   \nonumber \\ & & +\,C(0)\,t\int_0^t\ln P_\infty(\tau)\,P_\infty(\tau)\,d\tau
   \nonumber \\ & \le & C(0)P_\infty(t)^{1/2}\int_0^t
   \,P_\infty(\tau)^{1/2}\,m_2(\tau)^{1/2}\,d\tau
   \nonumber \\ & & +\,C(0)\,(1+t)\int_0^t\ln P_\infty(\tau)\,P_\infty(\tau)\,d\tau.
\end{eqnarray}
Next, by the definition of $E_\sharp$,
\begin{eqnarray*}
   I_\sharp(t) & = & \int_0^t V_0(s)\cdot E_\sharp(s, X_0(s))\,ds
   \\ & = & -\int_0^t ds\int_{|y|\le s}\frac{dy}{|y|}
   \int_{\R^3} dp\,V_0(s)\cdot K_{E,\,\sharp}(\omega, v)\,(Lf)(s-|y|, X_0(s)+y, p).
\end{eqnarray*}
From (\ref{anglr}) in Lemma \ref{kern-bd} it hence follows that
\begin{eqnarray*}
   |I_\sharp(t)| & \le & \int_0^t ds\int_{|y|\le s}\frac{dy}{|y|}
   \int_{\R^3}\frac{dp}{\sqrt{1+p^2}}\,\frac{1}{1+v\cdot\omega}
   \,(|E|+|B|)(s-|y|, X_0(s)+y)\,
   \\ & & \hspace{9em} \times f(s-|y|, X_0(s)+y, p)
   \\ & \le & \int_0^t ds\int_{|y|\le s}\frac{dy}{|y|}
   \,(|E|+|B|)(s-|y|, X_0(s)+y)\,\sigma_{-1}(s-|y|, X_0(s)+y)
   \\ & = & \int_0^t d\tau\int_{\tau}^t ds\,(s-\tau)\int_{|\omega|=1} dS(\omega)
   \,(|E|+|B|)(\tau, X_0(s)+(s-\tau)\omega)\,
   \\ & & \hspace{12em} \times\sigma_{-1}(\tau, X_0(s)+(s-\tau)\omega).
\end{eqnarray*}
Next the transformation $\Phi_\tau$ from (\ref{palla-trans}) is used
on $M_\tau=[\tau, t]\times [0, 2\pi] \times [0, \pi]$. This yields
\begin{eqnarray*}
   |I_\sharp(t)| & \le & \int_0^t d\tau\int_{\tau}^t ds\,(s-\tau)\int_0^{2\pi} d\varphi
   \int_0^\pi d\theta\,\sin\theta
   \,(|E|+|B|)(\tau, X_0(s)+(s-\tau)\omega)\,
   \\ & & \hspace{14em} \times\sigma_{-1}(\tau, X_0(s)+(s-\tau)\omega)
   \\ & = & \int_0^t d\tau\int_{\Phi_\tau(M_\tau)}dy\,\frac{1}{(s-\tau)}
   \,\frac{1}{(1+V_0(s)\cdot\omega)}
   \,(|E|+|B|)(\tau, y)\,\sigma_{-1}(\tau, y)
\end{eqnarray*}
for $s=s(y)$ and $\omega=\omega(y)$. Now fix $\eps\in ]0, \frac{1}{20}]$
and define $\alpha_\eps=\frac{12(1-\eps)}{7-20\eps+4\eps^2}\in ]1, 2[$.
The general H\"older inequality in $y$ for the exponents
$(\alpha_\eps, \frac{6}{1+2\eps}, \frac{4(1-\eps)}{1+2\eps})$
in conjunction with (\ref{palla-trans}) and Lemma \ref{efu-lem1}(c)
implies that
\begin{eqnarray*}
   |I_\sharp(t)| & \le & \int_0^t d\tau\,
   \bigg(\int_{\Phi_\tau(M_\tau)} dy\,\frac{1}{{(s-\tau)}^{\alpha_\eps}}
   \,\frac{1}{{(1+V_0(s)\cdot\omega)}^{\alpha_\eps}}\bigg)^{1/\alpha_\eps}
   \\ & & \hspace{4em} \times {\|(|E|+|B|)(\tau)\|}_{L^{\frac{6}{1+2\eps}}_x(\R^3)}
   \,{\|\sigma_{-1}(\tau)\|}_{L^{\frac{4(1-\eps)}{1+2\eps}}_x(\R^3)}
   \\ & = & \int_0^t d\tau\,
   \bigg(\int_{\tau}^t ds\int_0^{2\pi} d\varphi
   \int_0^\pi d\theta\,\sin\theta
   \,\frac{(s-\tau)^{2-\alpha_\eps}}{(1+V_0(s)\cdot\omega)^{\alpha_\eps-1}}\bigg)^{1/\alpha_\eps}
   \\ & & \hspace{4em} \times {\|(|E|+|B|)(\tau)\|}_{L^{\frac{6}{1+2\eps}}_x(\R^3)}
   \,{\|\sigma_{-1}(\tau)\|}_{L^{\frac{4(1-\eps)}{1+2\eps}}_x(\R^3)}
   \\ & \le & t^{\frac{2}{\alpha_\eps}-1}\int_0^t d\tau\,
   \bigg(\int_{\tau}^t ds\int_{|\omega|=1}
   \,\frac{dS(\omega)}{(1+V_0(s)\cdot\omega)^{\alpha_\eps-1}}\bigg)^{1/\alpha_\eps}
   \\ & & \hspace{4em} \times {\|(|E|+|B|)(\tau)\|}_{L^{\frac{6}{1+2\eps}}_x(\R^3)}
   \,{\|\sigma_{-1}(\tau)\|}_{L^{\frac{4(1-\eps)}{1+2\eps}}_x(\R^3)}
   \\ & \le & C(\eps)\,t^{\frac{3}{\alpha_\eps}-1}\int_0^t
   {\|(|E|+|B|)(\tau)\|}_{L^{\frac{6}{1+2\eps}}_x(\R^3)}
   \,{\|\sigma_{-1}(\tau)\|}_{L^{\frac{4(1-\eps)}{1+2\eps}}_x(\R^3)}\,d\tau.
\end{eqnarray*}
Returning to (\ref{P0bd}), it follows from this estimate and (\ref{Iflat-bd}) that
\begin{eqnarray*}
   {\bf 1}_{\{tP_\infty(t)\ge 1\}}\,|P_0(t)|
   & \le & \sqrt{1+P_0(0)^2}+\int_0^t V_0(s)\cdot (E_D+E_{DT})(s, X_0(s))\,ds
   \\ & & +\,{\bf 1}_{\{tP_\infty(t)\ge 1\}}\,|I_\flat(t)|+{\bf 1}_{\{tP_\infty(t)\ge 1\}}\,|I_\sharp(t)|
   \\ & \le & \sqrt{1+P_0(0)^2}+({\rm data})+C(0)P_\infty(t)^{1/2}\int_0^t
   \,P_\infty(\tau)^{1/2}\,m_2(\tau)^{1/2}\,d\tau
   \\ & & +\,C(0)\,(1+t)\int_0^t\ln P_\infty(\tau)\,P_\infty(\tau)\,d\tau
   \\ & & +\,C(\eps)\,t^{\frac{3}{\alpha_\eps}-1}\int_0^t
   {\|(|E|+|B|)(\tau)\|}_{L^{\frac{6}{1+2\eps}}_x(\R^3)}
   \,{\|\sigma_{-1}(\tau)\|}_{L^{\frac{4(1-\eps)}{1+2\eps}}_x(\R^3)}\,d\tau.
\end{eqnarray*}
Since $(X_0, P_0)$ is in the support of $f$, $\sqrt{1+P_0(0)^2}\le ({\rm data})$
uniformly in the characteristic, as $f^{(0)}\in C_0^1(\R^3\times\R^3)$
is compactly supported in $p$ by assumption; also see (\ref{charinsupp}).
If $f(s, x, p)\neq 0$, then $f(s, x, p)=f(s, X_0(s), P_0(s))$
for a characteristic $(X_0, P_0)$ in the support of $f$.
It follows from the preceding estimate that
\begin{eqnarray}\label{herma}
   {\bf 1}_{\{tP_\infty(t)\ge 1\}}\,P_\infty(t)
   & \le & ({\rm data})+C(0)P_\infty(t)^{1/2}\int_0^t
   \,P_\infty(\tau)^{1/2}\,m_2(\tau)^{1/2}\,d\tau
   \nonumber \\ & & +\,C(0)\,(1+t)\int_0^t\ln P_\infty(\tau)\,P_\infty(\tau)\,d\tau
   \nonumber \\ & & +\,C(\eps)\,t^{\frac{3}{\alpha_\eps}-1}\int_0^t
   {\|(|E|+|B|)(\tau)\|}_{L^{\frac{6}{1+2\eps}}_x(\R^3)}
   \,{\|\sigma_{-1}(\tau)\|}_{L^{\frac{4(1-\eps)}{1+2\eps}}_x(\R^3)}\,d\tau
   \qquad
\end{eqnarray}
for $t\in [0, T_{\rm max}[$. Now suppose that $T_{\rm max}<\infty$.
Then by the assumption (\ref{criter-hyp}),
\[ {\|\sigma_{-1}\|}_{L^\infty_t L^2_x(S_T)}\le\varpi(T)
   \le\max_{T'\in [0, T_{{\rm max}}]}\varpi(T')=:\varpi_{{\rm max}}<\infty \]
for all $T\in [0, T_{{\rm max}}[$. Thus according to Corollaries \ref{E+B},
\ref{sig-1L4}, and \ref{mom-cor} for $T, t\in [0, T_{{\rm max}}[$,
\begin{eqnarray}
   {\Big\||E|+|B|\Big\|}_{L^\frac{3}{1-\eps}_t L^\frac{6}{1+2\eps}_x(S_T)}
   & \le & C_6(\eps, T, {\rm data}, {\|\sigma_{-1}\|}_{L^\infty_t L^2_x(S_T)})
   \nonumber \\ & \le & C_6(\eps, T_{{\rm max}}, {\rm data}, \varpi_{{\rm max}}),
   \label{repuar} \\[1ex] {\|\sigma_{-1}(t)\|}_{L^{\frac{4(1-\eps)}{1+2\eps}}_x(\R^3)}
   & \le & C_8(0, \eps, t, {\rm data}, {\|\sigma_{-1}\|}_{L^\infty_t L^2_x(S_t)})
   \,P_\infty(t)
   \nonumber \\ & \le & C_8(0, \eps, T_{{\rm max}}, {\rm data}, \varpi_{{\rm max}})
   \,P_\infty(t),
   \nonumber \\[1ex] m_2(t) & \le & m_{\frac{3(1-2\eps)}{1+2\eps}}(t)
   \le C_7(0, \eps, t, {\rm data}, {\|\sigma_{-1}\|}_{L^\infty_t L^2_x(S_t)})
   \nonumber \\ & \le & C_7(0, \eps, T_{{\rm max}}, {\rm data}, \varpi_{{\rm max}}).
   \nonumber
\end{eqnarray}
Henceforth the dependence of the constants on ${\cal L}(0)$, the fixed $\eps$,
and the initial data is suppressed, and only the dependence on $T_{{\rm max}}$
and $\varpi_{{\rm max}}$ is made explicit. Thus (\ref{herma}) leads to
\begin{eqnarray*}
   {\bf 1}_{\{tP_\infty(t)\ge 1\}}\,P_\infty(t)
   & \le & ({\rm data})+C(T_{{\rm max}}, \varpi_{{\rm max}})
   \,P_\infty(t)^{1/2}\int_0^t\,P_\infty(\tau)^{1/2}\,d\tau
   \\ & & +\,C(0)\,(1+T_{{\rm max}})\int_0^t\ln P_\infty(\tau)\,P_\infty(\tau)\,d\tau
   \\ & & +\,C(T_{{\rm max}}, \varpi_{{\rm max}})\,T_{{\rm max}}^{\frac{3}{\alpha_\eps}-1}
   \int_0^t {\|(|E|+|B|)(\tau)\|}_{L^{\frac{6}{1+2\eps}}_x(\R^3)}\,P_\infty(\tau)\,d\tau
\end{eqnarray*}
for $t\in [0, T_{\rm max}[$. Since $P_\infty(\tau)\le P_\infty(t)^{1/2}P_\infty(\tau)^{1/2}$,
it follows that for a certain constant $C_2=C_2(T_{{\rm max}}, \varpi_{{\rm max}})>0$,
\begin{eqnarray*}
   \lefteqn{{\bf 1}_{\{tP_\infty(t)\ge 1\}}\,P_\infty(t)^{1/2}}
   \\ & \le & ({\rm data})+C(T_{{\rm max}}, \varpi_{{\rm max}})
   \,\int_0^t\,P_\infty(\tau)^{1/2}\,d\tau
   \\ & & +\,C(T_{{\rm max}}, \varpi_{{\rm max}})
   \int_0^t\ln P_\infty(\tau)^{1/2}\,P_\infty(\tau)^{1/2}\,d\tau
   \\ & & +\,C(T_{{\rm max}}, \varpi_{{\rm max}})
   \int_0^t {\|(|E|+|B|)(\tau)\|}_{L^{\frac{6}{1+2\eps}}_x(\R^3)}\,P_\infty(\tau)^{1/2}\,d\tau
   \\ & \le & ({\rm data})+C_2(T_{{\rm max}}, \varpi_{{\rm max}})
   \int_0^t\Big(1+{\|(|E|+|B|)(\tau)\|}_{L^{\frac{6}{1+2\eps}}_x(\R^3)}\Big)\,\ln P_\infty(\tau)^{1/2}
   \,P_\infty(\tau)^{1/2}\,d\tau
\end{eqnarray*}
for $t\in [0, T_{\rm max}[$. By the local existence theorem (Theorem \ref{locexi}),
there is a constant $C_1>0$ such that $\max_{t\in [0,\,T_{{\rm max}}/2]} P_\infty(t)
=P_\infty(T_{{\rm max}}/2)\le C_1$. Hence if $tP_\infty(t)\le 1$ and $t\in [0, T_{{\rm max}}/2]$, then
\[ P_\infty(t)^{1/2}\le C_1^{1/2}. \]
On the other hand, if $tP_\infty(t)\le 1$ and $t\in [T_{{\rm max}}/2, T_{{\rm max}}[$, then
\[ P_\infty(t)^{1/2}\le\frac{1}{\sqrt{t}}\le\bigg(\frac{2}{T_{{\rm max}}}\bigg)^{1/2}. \]
Therefore
\begin{eqnarray*}
   P_\infty(t)^{1/2}
   & \le & C_1(T_{{\rm max}}, \varpi_{{\rm max}})
   \\ & & +\,C_2(T_{{\rm max}}, \varpi_{{\rm max}})
   \int_0^t\Big(1+{\|(|E|+|B|)(\tau)\|}_{L^{\frac{6}{1+2\eps}}_x(\R^3)}\Big)\,\ln P_\infty(\tau)^{1/2}
   \,P_\infty(\tau)^{1/2}\,d\tau
\end{eqnarray*}
for $t\in [0, T_{\rm max}[$, where $C_1(T_{{\rm max}}, \varpi_{{\rm max}})
=({\rm data})+C_1^{1/2}+(\frac{2}{T_{{\rm max}}})^{1/2}$. This integral inequality
and (\ref{repuar}) imply that
\begin{eqnarray*}
   \ln P_\infty(t)^{1/2} & \le & C(T_{{\rm max}}, \varpi_{{\rm max}})
   \,\exp\bigg(\int_0^t\Big(1+{\|(|E|+|B|)(\tau)\|}_{L^{\frac{6}{1+2\eps}}_x(\R^3)}\Big)\,d\tau\bigg)
   \\ & \le & C(T_{{\rm max}}, \varpi_{{\rm max}})
   \,\exp\bigg(t+t^{\frac{2+\eps}{3}}
   \,{\Big\||E|+|B|\Big\|}_{L^\frac{3}{1-\eps}_t L^\frac{6}{1+2\eps}_x(S_t)}\bigg)
   \\ & \le & C(T_{{\rm max}}, \varpi_{{\rm max}})
   \,\exp\bigg(T_{{\rm max}}+T_{{\rm max}}^{\frac{2+\eps}{3}}
   \,C_6(\eps, T_{{\rm max}}, {\rm data}, \varpi_{{\rm max}})\bigg)
   \\ & \le & C_3(T_{{\rm max}}, \varpi_{{\rm max}})
\end{eqnarray*}
for $t\in [0, T_{\rm max}[$. Defining $\varpi_1(t)=\exp(C_3(T_{{\rm max}}, \varpi_{{\rm max}}))^2$,
the criterion (\ref{Pinfty-crit}) in Theorem \ref{locexi} is verified for $\varpi_1$.
From this result it hence follows that $T_{{\rm max}}=\infty$, which is a contradiction
to what was supposed before. As a consequence, $T_{{\rm max}}=\infty$ must be satisfied 
and the proof of Theorem \ref{sig-1thm} is complete. {\hfill$\Box$}\bigskip


\setcounter{equation}{0}

\section{Proof of Corollary \ref{thmcor}}
\label{prcor-sect} 

\begin{lemma} Define $\sigma_{-1}$ by (\ref{sig-1-def}) 
and ${\cal I}_\theta$ by (\ref{calidef}). 
Then for every $a\in [0, \infty[$ and $\eps>0$ there is a constant
$C=C(0, a, \eps)>0$ such that
\begin{equation}\label{knsasil}
   \sigma_{-1}(t, x)\le C\Big(1+{\cal I}_{a+1}(t, x)^{\frac{2+\eps a}{2+a}}\Big).
\end{equation}
\end{lemma}
{\bf Proof\,:} Fix $\omega\in\R^3$ such that $|\omega|=1$.
Since $1+v\cdot\omega\ge 1-|v|\ge\frac{1}{2(1+p^2)}$,
it follows for $R\in [10, \infty[$ and $\eps\in ]0, 2]$ that
\begin{eqnarray*}
   \lefteqn{\int_{\R^3}\frac{dp}{\sqrt{1+p^2}}
   \,\frac{1}{(1+v\cdot\omega)}\,f}
   \\ & = & \int_{|p|\le R,\,1+v\cdot\omega\le\eps}\frac{dp}{\sqrt{1+p^2}}
   \,\frac{1}{(1+v\cdot\omega)}\,f
   +\int_{|p|\le R,\,1+v\cdot\omega>\eps}\frac{dp}{\sqrt{1+p^2}}
   \,\frac{1}{(1+v\cdot\omega)}\,f
   \\ & & +\,\int_{|p|>R}\frac{dp}{\sqrt{1+p^2}}
   \,\frac{1}{(1+v\cdot\omega)}\,f
   \\ & \le & 2{\cal L}(0)\int_{|p|\le R,\,1+v\cdot\omega\le\eps} dp\,\sqrt{1+p^2}
   +{\cal L}(0)\int_{|p|\le R,\,1+v\cdot\omega>\eps}\frac{dp}{\sqrt{1+p^2}}
   \,\frac{1}{(1+v\cdot\omega)}
   \\ & & +\,2\int_{|p|>R} dp\,\sqrt{1+p^2}\,f
\end{eqnarray*}
For the first integral, the transformation (\ref{sigma-r}) yields
\begin{eqnarray*}
   \int_{|p|\le R,\,1+v\cdot\omega\le\eps} dp\,\sqrt{1+p^2}
   & = & \int_{|p|\le R,\,1+v_3\le\eps} dp\,\sqrt{1+p^2}
   \\ & \le & C\int_0^R dr\,r^2\sqrt{1+r^2}\int_0^{\pi} d\theta\,\sin\theta
   \,{\bf 1}_{\{1+\frac{r\cos\theta}{\sqrt{1+r^2}}\le\eps\}}
   \\ & \le & C\int_0^{R^\flat}\frac{d\sigma\,\sigma^2}{(1-\sigma^2)^3}
   \,\int_{-1}^1 ds\,{\bf 1}_{\{1+\sigma s\le\eps\}}
   \\ & \le & C\eps\int_0^{R^\flat}\frac{d\sigma\,\sigma}{(1-\sigma^2)^3}\le C\eps R^4.
\end{eqnarray*}
Similarly, the second integral can be bounded by
\begin{eqnarray*}
   \int_{|p|\le R,\,1+v\cdot\omega>\eps}\frac{dp}{\sqrt{1+p^2}}
   \,\frac{1}{(1+v\cdot\omega)}
   & \le & C\int_0^{R^\flat}\frac{d\sigma\,\sigma^2}{(1-\sigma^2)^2}
   \,\int_{-1}^1\frac{ds}{1+\sigma s}\,{\bf 1}_{\{1+\sigma s>\eps\}}
   \\ & \le & C\int_0^{R^\flat}\frac{d\sigma\,\sigma}{(1-\sigma^2)^2}
   \,\int_{\eps}^2\frac{d\tau}{\tau}\le C\ln\Big(\frac{2}{\eps}\Big)R^2.
\end{eqnarray*}
Hence for $a\in [0, \infty[$,
\[ \sigma_{-1}\le C(0)\eps R^4+C(0)\ln\Big(\frac{2}{\eps}\Big)R^2
   +2R^{-a}\int_{|p|>R} (1+p^2)^{\frac{a+1}{2}}\,f\,dp. \]
Upon choosing $\eps=1/R^2$, this leads to
\begin{equation}\label{pcikkleb}
   \sigma_{-1}\le C(0)\ln(R)R^2
   +2R^{-a}\int_{\R^3} (1+p^2)^{\frac{a+1}{2}}\,f\,dp.
\end{equation}
Furthermore, $\sigma_{-1}\le 2\int_{\R^3}\sqrt{1+p^2}\,f\,dp$
is always satisfied. Next fix a constant $C_\ast=C_\ast(a)$ such that
\[ I^{\frac{1}{2+a}}(\ln(10+I))^{-\frac{1}{2+a}}\ge 10,
   \quad I\ge C_\ast(a). \]
If ${\cal I}_{a+1}\le C_\ast(a)$, then
\[ \sigma_{-1}\le 2\int_{\R^3}\sqrt{1+p^2}\,f\,dp\le 2{\cal I}_{a+1}\le 2C_\ast(a). \]
On the other hand, if ${\cal I}_{a+1}\ge C_\ast(a)$,
then take $R={\cal I}_{a+1}^{\frac{1}{2+a}}(\ln(10+{\cal I}_{a+1}))^{-\frac{1}{2+a}}\ge 10$
in (\ref{pcikkleb}) to obtain
\begin{eqnarray*}
   \sigma_{-1} & \le & C(0)\ln(R)R^2+2R^{-a}{\cal I}_{a+1}
   \\ & \le & C(0, a)\,{\cal I}_{a+1}^{\frac{2}{2+a}}\Big([\ln {\cal I}_{a+1}
   -\ln\ln (10+{\cal I}_{a+1})]
   \,(\ln(10+{\cal I}_{a+1}))^{-\frac{2}{2+a}}+(\ln(10+{\cal I}_{a+1}))^{\frac{a}{2+a}}\Big)
   \\ & \le & C(0, a)\,{\cal I}_{a+1}^{\frac{2}{2+a}}\,(\ln(10+{\cal I}_{a+1}))^{\frac{a}{2+a}}.
\end{eqnarray*}
If $\eps>0$ is fixed, then select $C=C(a, \eps)$ such that $\ln(10+I)\le C(a, \eps)I^\eps$
whenever $I\ge C_\ast(a)$. Therefore ${\cal I}_{a+1}\ge C_\ast(a)$ yields
\[ \sigma_{-1}\le C(0, a)\,C(a, \eps)^{\frac{a}{2+a}}
   \,{\cal I}_{a+1}^{\frac{2+\eps a}{2+a}}, \]
and hence (\ref{knsasil}). {\hfill$\Box$}\bigskip

\begin{cor}\label{mudase} 
For every $a\in [0, \infty[$ and $\eps>0$ there is a constant $C=C(0, a, \eps)>0$ such that
\[ {\|\sigma_{-1}(t)\|}_{L^2_x(\R^3)}^2\le C\Big(1+t^3+{\|{\cal I}_{a+1}(t)\|}_{L^q_x(\R^3)}^q\Big) \] 
for $q=\frac{2(2+\eps a)}{2+a}$. 
\end{cor}
{\bf Proof\,:} Let $R_0$ be fixed such that $f^{(0)}(x, p)=0$ for $|x|\ge R_0$. 
Then (\ref{charg}) and (\ref{ff0}) implies that $f(t, x, p)=0$ for $|x|\ge R_0+t$. 
In particular, $\sigma_{-1}(t, x)=0$ for $|x|\ge R_0+t$. Hence squaring and integrating (\ref{knsasil}) 
it follows that 
\[ \int_{\R^3}\sigma_{-1}(t, x)^2\,dx\le C\Big((R_0+t)^3+\int_{\R^3} {\cal I}_{a+1}(t, x)^{\frac{2(2+\eps a)}{2+a}}\,dx\Big), \] 
which yields the claim. {\hfill$\Box$}\bigskip 

\noindent 
{\bf Proof of Corollary \ref{thmcor}\,:} Let $\theta>1$ and $q\in ]\frac{4}{\theta+1}, \infty[$ be given, 
and suppose that there is a function $\varpi_7\in C([0, \infty[)$ such that 
${\|{\cal I}_\theta(t)\|}_{L_x^q(\R^3)}\le\varpi_7(t)$ is verified for $t\in [0, T_{{\rm max}}[$. Defining 
\[ a=\theta-1,\quad\eps=\frac{q(2+a)}{2a}-\frac{2}{a}, \] 
we have $a>0$ and $\eps>0$, and in addition $q=\frac{2(2+\eps a)}{2+a}$. 
Thus we can apply Corollary \ref{mudase} to deduce that for $t\in [0, T_{{\rm max}}[$: 
\[ {\|\sigma_{-1}(t)\|}_{L^2_x(\R^3)}^2\le C\Big(1+t^3+{\|{\cal I}_\theta(t)\|}_{L^q_x(\R^3)}^q\Big)
   \le C\Big(1+t^3+\varpi_7(t)^q\Big). \]  
Hence Theorem \ref{sig-1thm} applies. {\hfill$\Box$}\bigskip 


\end{document}